\documentstyle[12pt]{article}
\begin{document}

\author{S.V. Ludkovsky}

\title{Structure of wrap groups of non-archimedean fibers.}

\date{17 November 2011 \thanks{Mathematics subject classification
(1991 Revision) 22-99, 46S10, 54H15 and 57S20.}\thanks{keywords:
wrap group; fiber bundle; infinite field; non-archimedean norm;
Cayley-Dickson algebra}} \maketitle

\begin{abstract}
This article is devoted to the investigation of structure of wrap
groups of fiber bundles over ultra-normed infinite fields and more
generally over Cayley-Dickson algebras. Iterated wrap groups are
studied as well. Their smashed products are constructed and studied.
\end{abstract}

\section{Introduction.} This article is a continuation of the previous one \cite{lulal12},
where wrap groups of non-archimedean fiber bundles were defined and
their existence was proved and some their properties were outlined.
Below other structure theorems are formulated and proved. \par  Wrap
groups of fiber bundles considered in this paper are constructed
with the help of families of mappings from a fiber bundle with
marked points into another fiber bundle with a marked point over an
ultra-normed field $\bf K$ or the octonion algebra ${\cal A}_3$
modeled on $\bf K$. This paper continues previous works  on this
theme, where generalized loop groups of manifolds over $\bf R$, $\bf
C$ and $\bf H$ were investigated, but neither for fibers nor over
octonions \cite{ludan,lugmlg,lujmslg,lufoclg}.

\par Applications of
quaternions in mathematics and physics can be found in
\cite{emch,guetze,lawmich}.
\par In this article wrap groups of different classes of smoothness are
considered. \par In particular, geometric loop groups have important
applications in modern physical theories (see \cite{ish,mensk} and
references therein). Groups of loops are also intensively used in
gauge theory. Wrap groups defined with the help of families of
mappings from a manifold $M$ into another manifold $N$ with a
dimension $dim_{\bf K} (M) >1$ over $\bf K$ can be used in the
membrane theory which is the generalization of the string
(superstring) theory.
\par In Section 2 smashed products of wrap groups are constructed.
Iterated wrap groups are studied as well. Their structure is
investigated in more details than in the previous paper. Relations
with path groups are studied. The main results of Section 2 are
Theorems 1, 5, 7, 8, 17, 18, 20, Propositions 6, 10, 11, 14, 16 and
Corollary 9.

\par All main results of this paper are obtained for the first time.

\section{Structure of wrap groups.}
\par Henceforth notations of the previous article are used \cite{lulal12}.

\par {\bf 1. Theorem.} {\it Let $M$ and $N$ be $C^{\alpha }_{\beta }$ manifolds
over an infinite field $\bf K$ or a Cayley-Dickson algebra ${\cal
A}_r$ for $\bf K$, $1\le r \le 3$, with a non-trivial multiplicative
ultra-norm, where $M$ and $N$ are of dimensions over $\bf K$ not
less than one, $dim_{\bf K} M\ge 1$ and $dim_{\bf K} N \ge 1$. Then
a wrap group $(W^MN)_{\alpha ,\beta }$ has no any nontrivial
continuous one parameter subgroup $\{ g^b: ~ b\in ({\bf K}, +) \}$.}
\par {\bf Proof.} Each manifold over ${\cal A}_r$ has a structure of
a manifold over $\bf K$ as well, hence it is sufficient to
demonstrate this theorem for $M$ and $N$ over $\bf K$. We shall
demonstrate that any non unit element, $g\ne e$, in $(W^MN)_{ \alpha
,\beta }$ does not belong to any one-parameter subgroup ${\bf K}\ni
b\mapsto g^b$, where $({\bf K},+)$ is considered as an additive
group. We already know that $(W^MN)_{\alpha ,\beta }$ is the
commutative group (see Theorem 6 \cite{lulal12}). As usually a
continuous one parameter subgroup means a continuous homomorphism
$\phi (b) = g^b$ so that $\phi (a+b)=\phi (a)+\phi (b)$ for all $a,
b \in {\bf K}$, $\phi : {\bf K}\to (W^MN)_{\alpha ,\beta }$ is
continuous, that is $g^{a+b}=g^a+g^b$, when the group operation in
the wrap group is denoted as the addition.
\par Suppose the contrary, that $\{ g^b: b\in ({\bf K},+) \} $ is
a nontrivial continuous one parameter subgroup, that is, $g^b\ne e$
for $b\ne 0$. The element $g$ is different from $e$, hence its
equivalence class $g = <f>_{K,\alpha ,\beta }$ for some $f\in
C^{\alpha }_{\beta }(({\bar M}, \{ s_{0,q}: q=1,...,k \}); (N,y_0))$
is different from $<w_0>_{K,\alpha ,\beta }$, where $w_0(M)= \{ y_0
\} $, $s_{0,q}$ with $q=1,...,k$ are marked points in $\bar M$ and
$y_0$ is a marked point in $N$, $M = {\bar M}\setminus \{ s_{0,q}:
q=1,..., k \} )$. In accordance with our convention the manifolds
$M$ and $N$ are modeled on locally $\bf K$ convex space $X$ and $Y$
over $\bf K$. Take a chart $(U_i^N, \mbox{}_N\phi _i)$ of $N$ so
that $y_0\in U_i^N\subset N$ (see \S 2.3 in \cite{lulal12}). Its
image $z_0 = \mbox{}_N\phi _i(y_0)$ belongs to $Y$.
 Without loss of generality we can consider that $z_0$
corresponds to zero in $Y$ making a shift $Y\ni y\mapsto y - z_0\in
Y$ in a case of necessity. Therefore, there exists a continuous
ultra-norm $v$ on $Y$ so that \par $(1)$  $t := \sup_{j\in \Lambda
_N, k\in \Lambda _M; x\in M_{j,k}} v( f_{j,k}(x) )
>0$, \\ where $f_{j,k}$ denotes a composition of $f$ with chart
mappings of the manifolds $N$ and $M$ respectively, $f_{j,k} =
\mbox{}_N\phi _j\circ f \circ \mbox{}_M\phi _k^{-1}$, $M_{j,k}$ is a
definition domain of $f_{j,k}$, $M_{j,k}\subset X$.
\par Let at first the field $\bf K$ be of zero characteristic $char
({\bf K})=0$. Then there exists a prime number $p>1$ so that
$0<|p|_{\bf K}<1$, consequently, $\lim_{n\to \infty }p^n =0$ in $\bf
K$, where $n\in \bf N$. If $g^b$ is continuous by $b$, then
$\lim_{n\to \infty } g^{p^n}=0$. On the other hand, \par $(2)$
$g^{p^n} = <\chi ^*(f^{\vee p^n})>_{K,\alpha ,\beta }$, \\  where
$f^{\vee m}$ denotes the $m$ fold wedge product of $f$ with itself,
$\chi ^*$ is the homomorphism for an $m$ fold wedge products of
mappings from $M^{\vee m}$ into $N$ induced by pairwise wedge
products by induction from Theorem 3.3 \cite{lulal12} denoted here
also by $\chi ^*$ \cite{lulal12}, $m\in \bf N$. But
\par $(3)$ $\sup_{j\in \Lambda _N, k\in \Lambda _M; x\in M_{j,k}} v(
[\chi^*(f^{\vee
p^n})]_{j,k}(x) ) = t >0$ \\
which contradicts $(2)$, consequently, the supposition about an
existence of a continuous nontrivial one parameter subgroup in
$(W^MN)_{\alpha ,\beta }$ was false.
\par If $\{ g^b: b\in ({\bf K},+) \} $ is a continuous non trivial one parameter subgroup and
$char ({\bf K})= z >0$, then the $z$ fold sum of the unit is zero,
$1+...+1 = 0 $, in the field $\bf K$. Thus we would have
$g^z=g^0=e$, since $g^0+g^b=g^b$ for each $b\in \bf K$. But again
\par $(4)$
$\sup_{j\in \Lambda _N, k\in \Lambda _M; x\in M} v( [\chi^*(f^{\vee
z})]_{j,k}(x) ) = t >0$ \\
that contradicts $(2)$. Therefore, the theorem is proved in this
case also.
\par {\bf 2. Remark.}  Mention that apart from the classical case
over the fields $\bf R$ or $\bf C$ the exponential function on
non-archimedean fields has only finite radius of convergence and
each differential equation with initial conditions generally have
infinite families of solutions, because for example on the field
$\bf Q_p$ of $p$-adic numbers there is an infinite family of
functions $f$ not equal to constants or even different from locally
constant functions, but with the derivative $df(x)/dx$ equal to zero
on $\bf Q_p$ \cite{sch1,sch2}.
\par It is also interesting to mention that under rather mild
conditions a manifold $M$ modeled on a non-archimedean Banach space
$X$ can be embedded into it as a clopen subset, $M\hookrightarrow X$
(see Reference $[24]$ in \cite{lulal12}).
\par {\bf 3. Iterated wrap groups.} We denote by $({\cal P}^ME)_{\alpha ,\beta }$ a space of equivalence
classes $<f>_{K,\alpha ,\beta }$ of $f\in C^{\alpha }_{\beta }({\bar
M},{\cal W})$ relative to the closures of orbits of the left action
of the family $Di^{\alpha }_{\beta ,0}(M)$ defined in \S 3.2
\cite{lulal12} of all generalized diffeomorphisms of a $C^{\alpha
'}_{\beta }$ differentiable space $\bar M$. This means, that $({\cal
P}^ME)_{\alpha ,\beta }$ is the quotient space of $C^{\alpha
}_{\beta }({\bar M},{\cal W})$ relative to the equivalence relation
$K_{\alpha ,\beta }$. It may be worthwhile to mention that in the
particular case of $\alpha =0$ the diffeomorphism group
$Diff^{\alpha }_{\beta ,0}({\bar M})$ reduces to the homeomorphism
group of preserving marked points mappings belonging to $C^0_{\beta
,0}$ class.
\par On the other hand, there is the embedding $\theta : C^{\alpha }_{\beta }
({\bar M},\{ s_{0,q}: q=1,...,k \} ; {\cal W},y_0) \hookrightarrow
C^{\alpha }_{\beta }({\bar M};{\cal W})$ and the evaluation mapping
${\hat ev}: C^{\alpha }_{\beta }({\bar M};{\cal W})\to N^k$ such
that ${\hat ev} (f) := ({\hat f}({\hat s}_{0,q}): q=k+1,..,2k)$,
${\hat ev}_{{\hat s}_{0,q}} (f) := {\hat f} ({\hat s}_{0,q})$, where
${\hat f} \in C^{\alpha }_{\beta }({\hat M};{\cal W})$ is such that
${\hat f} = f\circ \Xi $, $\Xi : {\hat M} \to \bar M$ is the
quotient mapping. Thus we get the diagram $C^{\alpha }_{\beta
}({\bar M},\{ s_{0,q}: q=1,...,k \} ; {\cal W},y_0)\to C^{\alpha
}_{\beta }({\bar M};{\cal W})\to N^k$ with $C^{\alpha }_{\beta }$
differentiable mappings, which induces the diagram
$\mbox{}_{l+1}C^{\alpha }_{\beta }({\bar M}, \{ s_{0,q}: q=1,...,k
\} ; {\cal W}, y_0)\to C^{\alpha }_{\beta }({\bar M},
\mbox{}_{l}C^{\alpha }_{\beta }({\bar M}, \{ s_{0,q}: q=1,...,k \} ;
{\cal W}, y_0)\to \mbox{}_{l}C^{\alpha }_{\beta }({\bar M}, \{
s_{0,q}: q=1,...,k \} ; {\cal W}, y_0)$ for each $l\in \bf N$, where
$\mbox{}_{l+1}C^{\alpha }_{\beta }({\bar M}, \{ s_{0,q}: q=1,...,k
\} ; {\cal W}, y_0) := C^{\alpha }_{\beta }({\bar M}, \{ s_{0,q}:
q=1,...,k \} ; \mbox{}_{l}C^{\alpha }_{\beta }({\bar M}, \{ s_{0,q}:
q=1,...,k \} ; {\cal W}, y_0))$, also $\mbox{}_{1}C^{\alpha }_{\beta
}({\bar M},\{ s_{0,q}: q=1,...,k \} ; {\cal W}, y_0) := C^{\alpha
}_{\beta }({\bar M},\{ s_{0,q}: q=1,...,k \} ; {\cal W}, y_0)$.
Using this procedure we get iterated wrap semigroups and groups
$(S^ME)_{l+1;\alpha ,\beta } := (S^M(S^ME)_{l;\alpha ,\beta
})_{\alpha ,\beta }$ and $(W^ME)_{l+1;\alpha ,\beta } :=
(W^M(W^ME)_{l;\alpha ,\beta })_{\alpha ,\beta }$, where
$(S^ME)_{1;\alpha ,\beta }:= (S^ME)_{\alpha ,\beta }$ and
$(W^ME)_{1;\alpha ,\beta } := (W^ME)_{\alpha ,\beta }$.
\par Evidently, if there are $C^{\alpha }_{\beta }$ and $C^{\alpha '}_{\beta }$
diffeomorphisms $\rho : {\bar M}\to {\bar M}_1$ and $\eta : N\to
N_1$ mapping marked points into respective marked points, then
$C^{\alpha }_{\beta }({\bar M},{\cal W})$ is isomorphic with
$C^{\alpha }_{\beta }({\bar M}_1,{\cal W}_1)$ and hence
$(W^ME)_{v;\alpha ,\beta }$ is $C^{\alpha }_{\beta }$ isomorphic as
the $C^{\alpha }_{\beta }$ differentiable space and $C^l_{\beta
}$-isomorphic as the $C^l_{\beta }$ differentiable Lie group with
$(W^{M_1}E_1)_{v;\alpha ,\beta }$ for each $v\in \bf N$, where
$l=\alpha ' - \alpha $, $l=\infty $ for $\alpha ' = \infty $,
$\alpha ' \ge \alpha \ge 0$. In particular, if $M$, $N$, $G$, $E$
are $C^{\alpha '}_{\beta }$ manifolds, then $(W^ME)_{\alpha ,\beta
}$ is a $C^{\alpha }_{\beta }$ manifold (see \S \S 3.3, 3.6
\cite{lulal12}). If $f: N\to N_1$ is a surjective map and $N_1$ is a
$C^{\alpha }_{\beta }$-differentiable space, then $N$ inherits a
structure of an $C^{\alpha }_{\beta }$-differentiable space with
plots having the local form $f\circ \rho : U\to N_1$, where $\rho :
U\to N$ is a plot of $N$.

\par {\bf 4. Lemma.} {\it Let $E$ be a $C^{\alpha '}_{\beta }$ principal bundle
and let $D$ be an everywhere dense subset in $N$ such that for each
$y\in D$ there exists an open neighborhood $V$ of $y$ in $N$ and a
differentiable map $p: V\to C^{\alpha }_{\beta }({\bar M}, \{
s_{0,q}: q=1,...,k \} ; V,y) := \{ f\in C^{\alpha }_{\beta }({\bar
M};V): f(s_{0,q}) =y, q=1,...,k \} $ such that ${\hat ev} _{{\hat
s}_{0,q}}({\hat p}(y))=y$ for each $q=1,...,2k$ and each $y\in N$,
where $p\circ \Xi = \pi \circ {\hat p}$. Then ${\hat ev}: C^{\alpha
}_{\beta }({\bar M};{\cal W})\to N^k$ is a $C^{\alpha }_{\beta }$
differentiable principal $(S^ME)_{\alpha ,\beta }$ bundle.}

\par {\bf Proof.} Let $ \{ (V_j, y_j): j\in J \} $ be a family such that
$y_j\in V_j\cap D$ for each $j$ and there exists a mapping $p_j:
V_j\to C^{\alpha }_{\beta }({\bar M}, \{ s_{0,q}: q=1,...,k \}
;V_j,y_j)$ so that ${\hat p}_j( {\hat s}_{0,q})(y) =y\times e$ for
each $q=1,...,2k$ and every $j$, where $ \{ V_j: j\in J \} $ is an
open covering of $N$, $y$ is a constant mapping from $\hat M$ into
$V_j$ with $y({\hat M})= \{ y \} $, where ${\hat p}_j({\hat
s}_{0,q})$ is the restriction to $V_j$ of the projection ${\hat
p}({\hat s}_{0,q}): ({\cal P}^ME)_{\alpha ,\beta }\to E$, while
$p_j(\Xi ({\hat x}))(y) = \pi \circ {\hat p}_j({\hat x})(y\times e)$
for each $y\in N$ and $x= \Xi ({\hat x})$ in $\bar M$, where ${\hat
x} \in \hat M$, $\Xi : {\hat M}\to \bar M$. Then $(W^ME)_{\alpha
,\beta }$ and $({\cal P}^ME)_{\alpha ,\beta }$ are supplied with the
$C^{\alpha }_{\beta }$-differentiable spaces structure (see Remark 3
above and Theorem 6 in \cite{lulal12}), where the embedding
$(S^ME)_{\alpha ,\beta }\hookrightarrow ({\cal P}^ME)_{\alpha ,\beta
}$ and the projection ${\hat ev}_{{\hat s}_{0,q}}: ({\cal
P}^ME)_{\alpha ,\beta }\to N$ are $C^{\alpha }_{\beta }$-maps.
\par Let a generalized diffeomorphism $\psi _j\in Di^{\alpha }_{\beta }(N)$ be
such that $\psi _j(y)=y_j$.
Specify a trivialization $\phi _j: {\hat p}_j^{-1}( {\hat s}_{0,q})
(V_j)\to V_j\times (S^ME)_{\alpha ,\beta }$ of the restriction
${\hat p}_j( {\hat s}_{0,q})|_{V_j}$ of the projection ${\hat p}_j(
{\hat s}_{0,q}): ({\cal P}^ME)_{\alpha ,\beta }\to E$ by the formula
$\phi _j(f) = (f({\hat s}_{0,q}), \psi _j\circ {\hat p}_j({\hat
s}_{0,q})(f))$ for each $f \in ({\cal P}^ME)_{\alpha ,\beta }$ with
$\pi \circ f ({\hat s}_{0,q})=y$, where $\psi _j\circ {\hat p}_j(f)
= \psi _j({\hat p}_j(f))$. Then $\phi _j^{-1}(y,g) = g^{-1}(\psi
_j\circ {\hat p}_j(y))=:\eta $, $\eta \in ({\cal P}^ME)_{\alpha
,\beta }$ with $\pi \circ \psi_j \circ f({\hat s}_{0,q})=y_j$, since
$G$ is a group, where $g= \psi _j\circ {\hat p}_j(f)$. A combining
of the family $\{ {\hat ev} _{{\hat s}_{0,q}}: q=k+1,...,2k \} $
induces a mapping ${\hat ev}: C^{\alpha }_{\beta }({\bar M}; {\cal
W})\to N^k$. By the construction above a fiber of this bundle is the
monoid $(S^ME)_{\alpha ,\beta }$.

\par {\bf 5. Theorem.} {\it If $\bar M$ and $N$ are $C^{\alpha }_{\beta }$
manifolds over ${\cal A}_r$, $0\le r \le 3$, $\hat M$ and $N$ are
embedded as clopen absolutely ${\cal A}_r$ convex subsets into
Banach spaces $Z_M$ and $Z_N$ over ${\cal A}_r$, then there exists a
$C^{\alpha }_{\beta }$-differentiable principal $(S^ME)_{\alpha
,\beta }$ bundle ${\hat ev}: ({\cal P}^ME)_{\alpha ,\beta }\to
N^k$.}
\par {\bf Proof.} In accordance with Lemma 4 it is sufficient to prove that
for each $y\in N$ there exist a neighborhood $U$ of $y$ in $N$ and a
$C^{\alpha }_{\beta }$-map $p_q: U\to C^{\alpha }_{\beta }({\bar
M},{\cal W})$ such that $ev_{s_{0,q}}(p_q(z))=z$ for each
$q=1,...,k$, $z\in U$, where $ev_x(f) = f(x)$.

\par  For $\hat M$ consider $\bf K$ linear mappings
$\zeta _q: B({\bf K},0,1)\to \hat M$ joining ${\hat s}_{0,q}$ with
${\hat s}_{0,q+k}$ such that $\zeta _q(0)= {\hat s}_{0,q}$ and
$\zeta _q(1)= {\hat s}_{0,q+k}$, since $\hat M$ is embedded as a
clopen absolutely ${\cal A}_r$ convex subset into a Banach space
$Z_M$ over ${\cal A}_r$, where $1\le q\le k$, $B(Z,z,R) := \{ y\in
Z: \rho (y,z)\le R \} $ denotes a ball of radius $0<R$ with a center
at $z$ in a metric space $Z$ with a metric $\rho $. We consider a
coordinate system $(x_1,...,x_m,...)$ in $\hat M$ over $\bf K$ such
that $x_1$ corresponds to a natural coordinate along $\zeta _q$.
This coordinate system is defined globally for a chosen $q$.
\par  The manifold $N$ is also embedded as the clopen absolutely ${\cal A}_r$ convex subset into
a Banach space $Z_N$ over ${\cal A}_r$. Then for each chart $U$ in
$N$ there exists a map $p_q: U\to ({\cal P}^MU)_{\alpha ,\beta }$
with $\pi \circ [p_q({\hat s}_{0,q+k})(z)] = z$ and $\pi \circ
[p_q({\hat s}_{0,q}(z)]=y$ for each $z\in U$, where $p_q\circ \zeta
_q=: {\hat \gamma }_{q,y,z}$ is a mapping so that $\pi \circ {\hat
\gamma }_{q,y,z}$ is a $\bf K$ linear mapping in $U$ joining $y$
with $z$, ${\hat \gamma }_{q,y,z}: B({\bf K},0,1)\to N$, ${\hat
\gamma }_{q,y,z}\circ \zeta ^{-1}_q(x_1)\in N$ for each $x_1$, where
$\pi : E\to N$ is the projection of the fiber bundle. Having
initially ${\hat \gamma }_{q,y,z}$ we extend it to ${\hat p}_q$ on
$\hat M$ with values in $E$ such that $p_q\circ \Xi = \pi \circ
{\hat p}_q$.

\par {\bf 6. Proposition.} {\it $(1)$.
The wrap group $(W^ME;N,G,{\bf P})_{\alpha ,\beta }$ has the a
structure of a principal $G^k$ bundle over $(W^MN)_{\alpha ,\beta }$
if either $M$ and $N$ satisfy conditions of Theorem 5 or $G^k$ acts
effectively on $(W^ME)_{\alpha ,\beta }$.
\par $(2)$. The abelianization $[(W^ME;N,G,{\bf P})_{\alpha ,\beta
}]_{ab}$ of the wrap group $(W^ME;N,G,{\bf P})_{\alpha ,\beta }$ is
isomorphic with $(W^ME;N,G_{ab},{\bf P})_{\alpha ,\beta }$.}

\par {\bf Proof. 1.} We have the bundle structure $\pi : E\to N$. It induces
the bundle structure ${\hat {\pi }}: (W^ME;N,G,{\bf P})_{\alpha
,\beta }\to (W^MN)_{\alpha ,\beta }$, since $\pi \circ {\bf
P}_{{\hat \gamma }, u} ={\hat \gamma }$. In accordance with Lemma 4
it is sufficient to show, that there exists a neighborhood $U_G$ of
$e$ in $(W^ME)_{\alpha ,\beta }$ and a $G$-equivariant mapping $\phi
: U_G\to (W^MN)_{\alpha ,\beta }$ (see Conditions 3.2$(P1-P4)$
\cite{lulal12}). Let $<{\bf P}_{{\hat {\gamma }},u}>_{\alpha ,\beta
}\in (W^ME)_{\alpha ,\beta }$, where ${\hat {\gamma }}: {\hat M}\to
N$, ${\hat {\gamma }} = \gamma \circ \Xi $, $\gamma : {\bar M}\to
N$, $\gamma (s_{0,q})=y_0$ for each $q=1,...,k$. Then $\pi \circ
{\bf P}_{{\hat {\gamma }},u}={\hat {\gamma }}$ and ${\bf P}_{{\hat
{\gamma }},u}$ is $G$-equivariant by the conditions defining the
parallel transport structure. This means that ${\bf P}_{{\hat
{\gamma }},u}(x)z = {\bf P}_{{\hat {\gamma }},uz}(x)$ for each $x\in
\hat M$ and $z\in G$ and every $u\in E_{y_0}$. We have that $uG =
\pi ^{-1}(y)$ for each $u\in E_y$ and $y\in N$.
\par Therefore, put $\phi = \pi _*$, where
$\pi _*<{\bf P}_{{\hat {\gamma }},u}>_{\alpha ,\beta } =<{{\hat
{\gamma }},u}>_{\alpha ,\beta }$ and take $U_G = \pi _*^{-1}(U)$,
where $U$ is a symmetric $U^{-1}=U$ neighborhood of $e$ in
$(W^MN)_{\alpha ,\beta }$.

\par The group $G$ acts effectively on $E$.
Then $G^k$ acts effectively on $(W^ME)_{\alpha ,\beta }$ if $N$ and
$M$ satisfy conditions of Theorem 5. Indeed, for each $\zeta _q$
from \S 5 there is $g_q\in G$ corresponding to ${\hat \gamma }({\hat
s}_{0,q+k})$ with ${\bf P}_{{\hat p}_q,{\hat s}_{0,q}\times e}({\hat
s}_{0,q+k}) =\{ y_0\times g_q \} \in E_{y_0}$, $g_q\in G$ for every
$q=1,...,k$. Moreover, $\pi _*^{-1}(\pi _*(<{\bf P}_{{\hat {\gamma
}},u}>_{\alpha ,\beta })) = <{\bf P}_{{\hat {\gamma }},u}>_{\alpha
,\beta }G^k $. Then the fibre of ${\hat \pi }: (W^ME;N,G,{\bf
P})_{\alpha ,\beta }\to (W^MN)_{\alpha ,\beta }$ is $G^k$. Due to
Conditions 2$(P1-P4)$ in Section 3.2 \cite{lulal12} it is the
principal $G^k$ differentiable bundle of class $C^{\alpha }_{\beta
}$.

\par {\bf 2.} Therefore, due to
$dim_{\bf K}N\ge 1$ the considered here wrap groups are infinite
dimensional over $\bf K$. Thus Statement $(2)$ follows from the
proof of Theorem 5 above, since the wrap groups $(W^MN)_{\alpha
,\beta }$ for $G=\{ e \} $ and $(W^ME;N,G_{ab},{\bf P})_{\alpha
,\beta }$ for $G=G_{ab}$ are commutative (see Theorem 6$(2)$
\cite{lulal12}).

\par {\bf 7. Theorem.} {\it Let $Diff^{\alpha '}_{\beta }(N)$ act transitively on $N$,
$\alpha \le \alpha '$, where $M$ and $N$ are embedded as clopen
${\cal A}_r$ absolutely convex subsets into Banach spaces $Z_M$ and
$Z_N$ over ${\cal A}_r$. For each $C^{\infty }$ manifold $N$ and an
$C^{\alpha }_{\beta }$ differentiable group $G$ such that ${\cal
A}_r^*\subset G$ with $1\le r\le 3$ there exists a homomorphism of
the $C^{\alpha }_{\beta }$ differentiable space of all equivalence
classes of $({\cal P}^ME)_{\alpha , \beta }$ relative to
$Diff^{\alpha '}_{\beta }(N)$ (see \S \S 1 and 2 in Section 3
\cite{lulal12} and \S 3 above) into $Hom^{\alpha }_{\beta
}((S^ME)_{\alpha , \beta },G^k)$. They are isomorphic, when $G$ is
commutative.}
\par {\bf Proof.} Mention that due to Theorem 5 the
$C^{\alpha }_{\beta }$-differentiable principal $(S^ME)_{\alpha ,
\beta }$ bundle ${\hat ev}: ({\cal P}^ME)_{\alpha , \beta }\to N^k$
has a parallel transport structure ${\hat {\bf P}}_{{\hat \gamma
},uz} (x)={\hat {\bf P}}_{{\hat {\gamma }},u}(x)z$ for each $x\in
\hat M$ and all $\gamma \in C^{\alpha }_{\beta }({\bar M},N)$ and
$u\in {\hat ev}^{-1}(\gamma (s_{0,k}))$ and every $z\in G$ and the
corresponding ${\hat \gamma }: {\hat M}\to N$ such that $\gamma
\circ \Xi = \hat {\gamma }$. If $x = {\hat s}_{0,q}$ with $1\le q\le
k$, then ${\hat {\bf P}}$ gives the identity homomorphism from
$(S^ME)_{\alpha , \beta }$ into $(S^ME)_{\alpha , \beta }$. If
$\theta : (S^ME)_{\alpha , \beta }\to G^k$ is an $C^{\alpha }_{\beta
}$ differentiable homomorphism, then the holonomy of the associated
parallel transport ${\hat {\bf P}}^{\theta }$ on the bundle $({\cal
P}^ME)_{\alpha , \beta }\times ^{\theta}G\to N^k$ is the
homomorphism $\theta : (S^ME)_{\alpha , \beta }\to G^k$ (see \S 6 in
Section 3 \cite{lulal12}). At the same time the group $G$ contains
continuous multiplicative one-parameter subgroups from ${\cal
A}_r^*$, where $1\le r\le 3$. If $g\in (W^MN)_{\alpha , \beta }$ and
$g\ne e$, then $g$ is of infinite order, since $w_0$ does not belong
to $g^n$ for each $n\ne 0$ non-zero integer $n$, where $w_0(M)= \{
y_0 \} $.

\par This holonomy induces a map $h: ({\cal P}^ME)_{\alpha , \beta }/{\cal Q} \to
Hom^{\alpha }_{\beta }((S^ME)_{\alpha , \beta },G^k)$ with values in
the family of homomorphisms of class $C^{\alpha }_{\beta }$ from
$(S^ME)_{\alpha , \beta }$ into $G^k$, where ${\cal Q}$ is an
equivalence relation caused by the transitive action of
$Diff^{\alpha '}_{\beta }(N)$ such that $(S^ME)_{\alpha , \beta }$
with distinct marked points $ \{ s_{0,q}: q=1,...,k \} $ in $M$ and
either $y_0$ or ${\tilde y}_0$ in $N$ are isomorphic, since there
exists $\psi \in Diff^{\alpha '}_{\beta }(N)$ such that $\psi
(y_0)={\tilde y}_0$.
\par If $G$ is commutative, then this map is the homomorphism, since
$(S^ME)_{\alpha , \beta }$ is the commutative monoid for a
commutative group $G$ (see Theorem 3 in Section 3 \cite{lulal12})
and $u{\bf P}_{{\hat {\gamma }}_1,v_1}(x_1) {\bf P}_{{\hat {\gamma
}}_2,v_2}(x_2)= u{\bf P}_{{\hat {\gamma }}_2,v_2}(x_2) {\bf
P}_{{\hat {\gamma }}_1,v_1}(x_1)$ for each $x_1, x_2\in \hat M$ and
$u, v_1, v_2 \in E_{y_0}$. We have the embedding $(S^ME)_{\alpha ,
\beta }\hookrightarrow (W^ME)_{\alpha , \beta }$. Thus a
homomorphism $\theta : (W^ME)_{\alpha , \beta }\to G^k$ has the
restriction on $(S^ME)_{\alpha , \beta }$ which is also the
homomorphism. \par For $G\supset {\cal A}_r^*$ there exists a family
of $f\in Hom^{\alpha }_{\beta }((S^ME)_{\alpha , \beta },G^k)$
separating elements of the wrap monoid $(S^ME)_{\alpha , \beta }$,
hence there exists the embedding of $(S^ME)_{\alpha , \beta }$ into
$Hom^{\alpha }_{\beta }((S^ME)_{\alpha , \beta },G^k)$. The bundle
$({\cal P}^ME)_{\alpha , \beta }\times ^{\theta }G\to N^k$ has the
induced parallel transport structure ${\bf P}^{\theta }$. The
holonomy of the parallel transport structure on $({\cal
P}^MN)_{\alpha , \beta }\times ^{\theta }G\to N^k$ is $\theta $.
Therefore, the map $C^{\alpha }_{\beta }((S^ME)_{\alpha , \beta
},G^k)\ni \theta \mapsto {\bf P}^{\theta }$ is inverse to $h$.

\par {\bf 8. Embeddings of wrap groups and normal subgroups.}
Suppose that \par $(E1)$ there are embeddings ${\bar
M}_2\hookrightarrow {\bar M}_1$ and ${\bar M} = {\bar M}_1\setminus
({\bar M}_2\setminus {\bar M}_{2,f})$ and ${\hat M}_2\hookrightarrow
{\hat M}_1$ and ${\hat M} = {\hat M}_1\setminus ({\hat M}_2\setminus
{\hat M}_{2,f})$ and $N_2\hookrightarrow N_1$ for $C^{\alpha
}_{\beta }$-manifolds with the same marked points $\{ s_{0,q}:
q=1,...,k \} $ for ${\bar M}_1$ and ${\bar M}_2$ and $\bar M$ and
$y_0\in N_2$ satisfying conditions of \S \S 1 and 2 in Section 3
\cite{lulal12} and \par $(E2)$ $G_2$ is a closed subgroup in $G_1$
with a complete relative to its uniformity principal fiber bundle
$E$ with a structure group $G_1$. Moreover, we suppose that atlases
of all embedded pairs of manifolds $A\hookrightarrow B$ are
consistent in the following sense. \par $(E3)$ Each chart $U$ of $A$
is contained in some chart $V$ of $B$ so that there exists a
$C^{\alpha }_{\beta }$ embedding $\theta _{U,V}: \phi
_U(U)\hookrightarrow \phi _V(V)\subset X_B$, where \par $(E4)$ $\phi
_U: U\to \phi _U(U)$ and $\phi _V: V\to \phi _V(V)$ are
homeomorphisms, \par $(E5)$ $\phi _U(U)$ and $\phi _V(V)$ are $\bf
K$ convex in $X_B$, where $X_B$ is a complete locally $\bf K$ convex
space. Moreover,
\par $(E6)$ there exists a topological Schauder basis $\{ e_j: j\in J \} $ in $X_B$,
where $J$ is a set. Suppose also that
\par $(E7)$ for each point $x_0\in \theta _{U,V}(\phi _U(U))$ either
$(x_0 + e_j{\bf K})\cap \theta _{U,V}(\phi _U(U))$ is clopen in
$\phi _V(V)\cap (x_0 + e_j{\bf K})$ or a singleton is. That is this
consistency of atlases is satisfied for pairs $({\bar M}_2, {\bar
M}_1)$, $({\hat M}_2, {\hat M}_1)$ and $(N_2,N_1)$. Suppose also
that \par $(E8)$ manifolds ${\hat M}_1$ and ${\bar M}_1$ are finite
dimensional over $\bf K$.

\par {\bf Theorem. 1}. {\it Then $(W^{M_2, \{ s_{0,q}: q=1,...,k \} }E;N_2,G_2,{\bf
P})_{\alpha , \beta }$ has an embedding as a closed subgroup into
$(W^{M_1, \{ s_{0,q}: q=1,...,k \} } E;N_1,G_1,{\bf P})_{\alpha ,
\beta }$.
\par {\bf 2}. The wrap group
$(W^{M_2, \{ s_{0,q}: q=1,...,k \} }E;N,G_2,{\bf P})_{\alpha , \beta
}$ is normal in \\ $(W^{M_1, \{ s_{0,q}: q=1,...,k \} } E;N,G_1,{\bf
P})_{\alpha , \beta }$ if and only if $G_2$ is a normal subgroup in
$G_1$.
\par {\bf 3}. In the latter case $(W^ME;N,G,{\bf P})_{\alpha , \beta }$
is isomorphic with \\ $(W^{M_1}E;N,G_1,{\bf P})_{\alpha , \beta
}/(W^{M_2}E;N,G_2,{\bf P})_{\alpha , \beta }$, where $G=G_1/G_2$.}

\par {\bf Proof. 1.} Manifolds ${\bar M}_1$ and ${\bar M}_2$ are
finite dimensional over the field $\bf K$. Each finite dimensional
vector space over the field $\bf K$ is isomorphic with ${\bf K}^n$
for some natural number $n\in \bf N$. Therefore, ${\bar M}_1$ and
${\bar M}_2$ have disjoint clopen atlases refining their initial
atlases. Without loss of generality we can take such atlases. \par
We recall that a system of $\bf K$ linearly independent vectors $ \{
e_j: j\in J \} $ is called a topological basis in a topological $\bf
K$ vector space $X$, if each $x\in X$ can be decomposed as a limit
of finite $\bf K$ linear combinations of elements $e_j$ with
components $a_j= \zeta _j(x)\in \bf K$ of $x$, where each $\zeta
_j(x)$ is a $\bf K$ linear functional on $X$, $x = \lim \sum_j \zeta
_j(x) e_j$. A topological basis $ \{ e_j: j\in J \} $ with
continuous $\bf K$ linear functionals $\zeta _j : X\to {\bf K}$ for
each $j\in J$ is called a Schauder basis of $X$.
\par For manifolds ${\hat M}_1$ and ${\bar M}_1$ a vector space $X$
on which they are modeled is finite dimensional over $\bf K$ by the
supposition of this theorem.
\par Moreover, $\phi _U(U)$ can be presented as a disjoint union $\phi
_U(U) = \bigcup_j A_j$ of balls clopen in $x_j + span_{\bf K} \{
e_{k_1(j)},...,e_{k_m(j)} \} $, where $m=n_2$ is a dimension of
$\phi _U(U)$ over $\bf K$, $k_1(j),...,k_m(j)$ are pairwise
different natural numbers in $ \{ 1,...,n_1 \} $, $n_1$ is a
dimension of $\phi _V(V)$ over $\bf K$, since $\phi _U(U)$ and $\phi
_V(V)$ satisfy Conditions $(E3-E8)$.
\par Each complete locally $\bf K$ convex space $Y$ is a projective
limit of Banach spaces $Y_{\xi }$ over $\bf K$, where $\xi \in
\Lambda _Y$, $\Lambda _Y$ is a directed set \cite{nari}. If $Y$ has
a topological Schauder basis, then each $Y_{\xi }$ has a topological
Schauder basis $ \{ e_{j,\xi }: j\in J_{\xi } \} $, where $J_{\xi }$
is a set. A function $f: {\hat M}_1\to Y$ is of $C^{\alpha }_{\beta
}$ class if and only if $\pi _{\xi }\circ f: {\hat M}_1\to Y_{\xi }$
is of $C^{\alpha }_{\beta }$ class for each $\xi $, where $\pi _{\xi
}: Y\to Y_{\xi }$ is the quotient mapping. Pointwise we have a
decomposition $\pi _{\xi }\circ f (x) = \lim \sum_j f_{j,\xi }(x)
e_{j,\xi }$, where $f_{j,\xi } (x) = \zeta _{j,\xi }(f(x))$, $x\in
{\hat M}_1$, $\zeta _{j,\xi } = \pi \circ \zeta _j$. If each
$f_{j,\xi }$ is of $C^{\alpha }_{\beta }$ class, then $\pi _{\xi
}\circ f$ is of $C^{\alpha }_{\beta }$ class.

\par If $f: T\to \bf K$ is a $C^{\alpha }_{\beta }$ function from a
clopen subset $T$ in ${\bf K}^m$ into $\bf K$ or $X_N$, then it has
a $C^{\alpha }_{\beta }$ extension on ${\bf K}^m$ taking $f|_{({\bf
K}^m\setminus T)} =g$, where $g$ is a $C^{\alpha }_{\beta }$
function on ${\bf K}^m\setminus T$, since the latter set is also
clopen in ${\bf K}^m$ (see \S 2 in Section 2 \cite{lulal12}). If $T$
is a singleton in $\bf K$, then evidently a locally constant or some
other $C^{\alpha }_{\beta }$ extension from $T$ on $\bf K$ exists.
\par Thus we get $C^{\alpha }_{\beta }$ extensions from $A_j$ on
$B_j$, where $B_j$ is clopen in $X_{M_1}$, $B_j\cap (x_j + span_{\bf
K} \{ e_{k_1(j)},...,e_{k_m(j)} \} ) = \theta _{U,V}(A_j)$ for each
$j$, $\bigcup_j B_j = \phi _V(V)$, due to Theorem 40
\cite{lujms157dif} for the finite dimensional over $\bf K$ space $X$
for functions with values in $X_{N_1}$, where $X_{N_1}$ is a $\bf K$
vector space on which $N_1$ is modeled. Combining these disjoint
clopen coverings and mappings on them we get a $C^{\alpha }_{\beta
}$ extension from $M_2$ onto $M_1$ . Therefore, if ${\hat {\gamma
}}_2\in C^{\alpha }_{\beta }({\hat M}_2, N_2)$, then it has a
$C^{\alpha }_{\beta }$ extension to ${\hat {\gamma }}_1\in C^{\alpha
}_{\beta }({\hat M}_1, N_1)$, since these manifolds are totally
disconnected and their atlases are consistent.
\par Thus the parallel transport structure ${\bf P}_{{\hat {\gamma
}}_1,u}$ over ${\hat M}_1$ serves as an extension of ${\bf P}_{{\hat
{\gamma }}_2,u}$ over ${\hat M}_2$. The uniform spaces $C^{\alpha
}_{\beta }({\bar M}_j, \{ s_{0,1},...,s_{0,k} \} ;{\cal W}_j,y_0)$
are complete for $j=1,2$, since the principal fiber bundle $E$ is
complete relative to its uniformity and the corresponding principal
fiber sub-bundle $E_2$ with the structure group $G_2$ is also
complete (see Theorem 8.3.6 \cite{eng}). Therefore, $C^{\alpha
}_{\beta }({\bar M}_2, \{ s_{0,1},...,s_{0,k} \} ; {\cal W}_2,y_0)$
has embedding as the closed subspace into $C^{\alpha }_{\beta
}({\bar M}_1, \{ s_{0,1},...,s_{0,k} \} ; {\cal W}_1,y_0)$. Using
Theorem 40 \cite{lujms157dif} as above we infer that each $C^{\alpha
}_{\beta }$ diffeomorphism of ${\bar M}_2$ has an $C^{\alpha
}_{\beta }$ extension to a diffeomorphism of ${\bar M}_1$. From the
condition that $G_2$ is a closed subgroup in $G_1$ we infer that
$(S^{M_2, \{ s_{0,q}: q=1,...,k \} }E;N_2,G_2,{\bf P})_{\alpha ,
\beta }$ has an embedding as a closed sub-monoid into $(S^{M_1, \{
s_{0,q}: q=1,...,k \} } E;N_1,G_1,{\bf P})_{\alpha , \beta }$ and
inevitably $(W^{M_2, \{ s_{0,q}: q=1,...,k \} }E;N_2,G_2,{\bf
P})_{\alpha , \beta }$ has an embedding as a closed subgroup into
$(W^{M_1, \{ s_{0,q}: q=1,...,k \} } E;N_1,G_1,{\bf P})_{\alpha ,
\beta }$ due to Theorem 6.1 in Section 3 \cite{lulal12}.

\par {\bf 2.} The groups $(W^{M_j, \{ s_{0,q}: q=1,...,k \} }N)_{\alpha , \beta }$
for $j=1, 2$ are commutative and $(W^{M_j, \{ s_{0,q}: q=1,...,k \}
}E)_{\alpha , \beta }$ is the $G_j^k$ principal fiber bundle on
$(W^{M_j, \{ s_{0,q}: q=1,...,k \} }N)_{\alpha , \beta }$ (see
Theorem 6.2 in Section 3 \cite{lulal12} and Proposition 6.1 above).
Therefore, $(W^{M_2, \{ s_{0,q}: q=1,...,k \} }E)_{\alpha , \beta }$
is the normal subgroup in $(W^{M_1, \{ s_{0,q}: q=1,...,k \}
}E)_{\alpha , \beta }$ if and only if $G_2$ is the normal subgroup
in $G_1$.

\par {\bf 3.} Consider the principal fiber bundle $E(N,G,\pi ,\Psi )$
with the structure group $G$ (see Section 2.6 \cite{lulal12}) and
the parallel transport structure ${\bf P}$ for the $C^{\alpha
}_{\beta }$ pseudo-manifold $\hat M$, where $G=G_1/G_2$ is the
quotient group. If ${\hat {\gamma }}_1\in C^{\alpha }_{\beta }({\hat
M}_1, N)$, then ${\hat {\gamma }}_1$ is the combination
\par $(i)$ ${\hat {\gamma }}_1 = {\hat {\gamma }}_2\nabla {\hat
{\gamma }}$,
\\ where ${\hat {\gamma }}_2$ and ${\hat {\gamma }}$ are
restrictions of ${\hat {\gamma }}_1$ on ${\hat M}_2$ and ${\hat M}$
respectively. We also have that each ${\hat {\gamma }}\in C^{\alpha
}_{\beta }({\hat M}, N)$ has an extension ${\hat {\gamma }}_1\in
C^{\alpha }_{\beta }({\hat M}_1, N)$. The manifold ${\hat M}_1$ is
metrizable by a metric $\rho $. For each $\epsilon
>0$ there exists $\psi \in Diff^{\alpha }_{\beta ,0}({\hat M}_1)$
such that $(\psi ({\hat M})\cap {\hat M}_2)\subset \bigcup_{l=1}^s
B({\hat M}_1,x_l,\epsilon )$ for some $x_l\in {\hat M}_1$ with
$l=1,...,s$ and $s\in \bf N$ and $\psi |_{{\hat M}_1\setminus ({\hat
M}\bigcup_{l=1}^s B({\hat M}_1,x_l,\epsilon _l))}=id$, since ${\hat
M}_1$ and ${\hat M}_2$ are finite dimensional manifolds, where
$B({\hat M}_1,x,R)$ denotes a ball in ${\hat M}_1$ containing a
point $x$ and of radius $0<R$, $0<\epsilon _l<\infty $. Therefore, a
using charts of the manifolds gives
\par $<{\bf P}_{{\hat {\gamma }},u}|_M>_{\alpha , \beta } = <{\bf
P}_{{\hat {\gamma }}_1,u}|_{M_1}>_{\alpha , \beta }/<{\bf P}_{{\hat
{\gamma }}_2,u}|_{M_2}>_{\alpha , \beta }$
\\ due to decomposition $(i)$, since ${\bf P}_{{\hat {\gamma
}},u}|_{M_j}\in G_j$ for $j=1,2$ and $G=G_1/G_2$ is the $C^{\alpha
'}_{\beta }$ quotient group with $\alpha '\ge \alpha $.
Consequently, $(W^ME;N,G,{\bf P})_{\alpha , \beta }$ is isomorphic with \\
$(W^{M_1}E;N,G_1,{\bf P})_{\alpha , \beta }/(W^{M_2}E;N,G_2,{\bf
P})_{\alpha , \beta }$ (see also \S \S 3, 6 in Section 3
\cite{lulal12}).

\par {\bf 9. Corollary.} {\it Let suppositions of Theorem 8 be
satisfied. Then the group $(W^MN)_{\alpha , \beta }$ is isomorphic
with the quotient group $(W^{M_1}N)_{\alpha , \beta }/
(W^{M_2}N)_{\alpha , \beta }$.}
\par {\bf Proof.} For $(W^MN)_{\alpha , \beta }$ taking $G=G_1=G_2=\{ e \} $
we get the statement of this corollary from Theorem 8.3.

\par {\bf 10. Proposition.} {\it  Suppose that ${\bar M} = {\bar M}_1\vee {\bar M}_2$, where
${\bar M}_1$ and ${\bar M}_2$ are $C^{\alpha }_{\beta }$
differentiable spaces satisfying Conditions 2.6, 3.1 and 3.2
\cite{lulal12} with the bunch taken by marked points $\{ s_{0,q}:
q=1,...,k \} $, then $(W^MN)_{\alpha , \beta }$ is isomorphic with
the internal direct product $(W^{M_1}N)_{\alpha , \beta }\otimes
(W^{M_2}N)_{\alpha , \beta }$.}

\par {\bf Proof.} The $C^{\alpha }_{\beta }$ differentiable space $\bar M$ has marked points
$\{ s_{0,q}: q=1,...,k \} $ such that $s_{0,q}$ corresponds to
$s_{0,q,1}$ glued with $s_{0,q,2}$ in the bunch ${\bar M}_1\vee
{\bar M}_2$ for each $q=1,...,k$, where $s_{0,q,j}\in {\bar M}_j$
are marked points $j=1, 2$. Each ${\bar M}_j$ satisfies Conditions
3.1$(S1-S5)$ and 3.2 in \cite{lulal12}, then $M$ satisfies them
also. Each $C^{\alpha ,w_0}_{\beta ,0}$ function on ${\bar M}_1$ in
$\bar M$ has a $C^{\alpha ,w_0}_{\beta ,0}$ extension as $w_0$ on
$M_2$. Due to the initial conditions at marked points $s_{0,q}$,
$q=1,...,k$, it has a $C^{\alpha , w_0}_{\beta ,0 }$ extension on
${\bar M}_2$ and thus on the entire $\bar M$ also. Therefore, quite
analogously to \S 8 $(W^{M_j, \{ s_{0,q}: q=1,...,k \} }N)_{\alpha ,
\beta }$ has an embedding as a closed subgroup into $(W^{M, \{
s_{0,q}: q=1,...,k \} } N)_{\alpha , \beta }$ for $j=1, 2$. If
$\gamma _j\in C^{\alpha }_{\beta }(({\bar M}_j, \{ s_{0,q}:
q=1,...,k \}); (N,y_0))$ for $j=1, 2$, then $\gamma _1\vee \gamma
_2\in C^{\alpha }_{\beta }(({\bar M}, \{ s_{0,q}: q=1,...,k \} );
(N,y_0))$. At the same time each $\gamma \in C^{\alpha }_{\beta
}(({\bar M}, \{ s_{0,q}: q=1,...,k \}); (N,y_0))$ has the
decomposition $\gamma = \gamma _1\vee \gamma _2$, where $\gamma _j =
\gamma |_{{\bar M}_j}$ for $j=1, 2$. Therefore, $<\gamma
>_{\alpha , \beta } = <\gamma _1\vee w_{0,2}>_{\alpha , \beta }\vee
<w_{0,1}\vee \gamma _2>_{\alpha , \beta }$, where $w_0(M)= \{ y_0 \}
$, $w_{0,j}= w_0|_{M_j}$ for $j=1, 2$, hence $(W^MN)_{\alpha , \beta
}$ is isomorphic with $(W^{M_1}N)_{\alpha , \beta }\otimes
(W^{M_2}N)_{\alpha , \beta }$.

\par {\bf 11. Propositions. 1.} {\it Let $\theta : N_1\to N$ be an
embedding with $\theta (y_1)=y_0$, or $F: E_1\to E$ be an embedding
of principal fiber bundles over ${\cal A}_r$ such that $\pi \circ
F|_{N_1\times e} =\theta \circ \pi _1$, then there exist embeddings
$\theta _*: (W^MN_1)_{\alpha , \beta }\to (W^MN)_{\alpha , \beta }$
and $F_*: (W^ME_1)_{\alpha , \beta }\to (W^ME)_{\alpha , \beta }$.}
\par {\bf 2.} {\it If $\theta : N_1\to N$ and $F: E_1\to E$ are a quotient
mapping and a quotient homomorphism such that $N_1$ is a covering
$C^{\alpha }_{\beta }$ differentiable space of a $C^{\alpha }_{\beta
}$ differentiable space $N$ satisfying conditions of \S 2.6
\cite{lulal12}, then $(W^MN)_{\alpha , \beta }$ is the quotient
group of some closed subgroup in $(W^MN_1)_{\alpha , \beta }$ and
$(W^ME)_{\alpha , \beta }$ is the quotient group of some closed
subgroup in $(W^ME_1)_{\alpha , \beta }$.}
\par {\bf 3.} {\it If there are a $C^{\alpha }_{\beta }$ diffeomorphism $f_1: M\to M_1$
and an $C^{\alpha '}_{\beta }$-isomorphism $f_2: E\to E_1$, then
wrap groups $(W^{M_1}E_1)_{\alpha , \beta }$ and $(W^ME)_{\alpha ,
\beta }$ are isomorphic.}

\par {\bf Proof. 1.} If $\gamma _1\in C^{\alpha }_{\beta }(({\bar M},\{ s_{0,q}: q=1,...,k \});
(N_1,y_1))$, then $\theta \circ \gamma _1=\gamma \in C^{\alpha
}_{\beta }(({\bar M}, \{ s_{0,q}: q=1,...,k \} ); (N,y_0))$,
$<\gamma
>_{\alpha , \beta } = \theta _*<\gamma _1>_{\alpha , \beta }$, where
$\theta _*<\gamma _1>_{\alpha , \beta } := \{ \theta \circ f:
fK_{\alpha , \beta }\gamma _1 \} $. In addition $F|_{E_{1,v}}$ gives
an embedding $F: G_1\to G$, where $G_1$ and $G$ are structural
groups of $E_1$ and $E$ correspondingly. Therefore, for the parallel
transport structures we get
\par $(1)$ $F\circ {\bf
P}^1_{{\hat {\gamma }}_1,v}(x) = {\bf P}_{{\hat {\gamma }},u}(x)$\\
for each $x\in \hat M$, where $F(v)=u$, $\pi \circ F = \theta \circ
\pi _1$, where ${\bf P}^1$ is for $E_1$ and $\bf P$ for $E$. Define
$F_*<{\bf P}^1_{{\hat {\gamma }}_1,v}>_{\alpha , \beta } := \{
F\circ g: gK_{\alpha , \beta }{\bf P}^1_{{\hat {\gamma }}_1,v} \} $.
Since $\theta $ and $F$ are $C^{\alpha }_{\beta }$ differentiable
mappings, then $\theta _*$ and $F_*$ are embeddings of $C^{\alpha
}_{\beta }$ differentiable spaces and group homomorphisms of
$C^l_{\beta }$ differentiable groups (see also Theorems 6 in Section
3 \cite{lulal12}).

\par {\bf 2.} By the conditions of this theorem $N_1$ is a covering of $N$,
that is each $y\in N$ has a neighborhood $V_y$ for which $\theta
^{-1}(V_y)$ is a disjoint union of open subsets in $N_1$. If an open
covering $\cal V$ of $\bar M$ and a function $f\in C^{\alpha
}_{\beta }({\bar M},N)$ are such that for each $\nu \in \cal V$
there exist $y\in N$ and $V_y$ as above for which the embedding
$f(\nu )\subset V_y$ is satisfied, then $f_1\in C^{\alpha }_{\beta
}({\bar M},N_1)$ exists so that $\theta \circ f_1 =f$. If $\gamma
\in C^{\alpha }_{\beta }(({\bar M}, \{ s_{0,q}: q=1,...,k \} );
(N,y_0))$, then there exists $\gamma _1 \in C^{\alpha }_{\beta
}(({\bar M}, \{ s_{0,q}: q=1,...,k \} ); (N_1,y_1))$ such that
$\theta \circ \gamma _1=\gamma $. This $\gamma _1$ evidently exists
due to total disconnectedness of $\bar M$ and $\gamma ({\bar M})$
and the choice axiom \cite{eng}, where $\gamma ({\bar M})\subset N$.
To each parallel transport in $E_1$ there corresponds a parallel
transport in $E$ so that Equation $(1)$ above is satisfied. Put
$\theta _*^{-1}<\gamma
>_{\alpha , \beta } = \{ <\gamma _1>_{\alpha , \beta }: \theta \circ
\gamma _1 = \gamma \} $ and $F_*^{-1}<{\bf P}_{{\hat {\gamma
}},u}>_{\alpha , \beta } := \{ <{\bf P}^1_{{\hat {\gamma
}}_1,v}>_{\alpha , \beta }: F\circ {\bf P}^1_{{\hat {\gamma }}_1,v}
= {\bf P}_{{\hat {\gamma }},u} \} $, where $F(v)=u$.
\par Thus we obtain quotient mappings $\theta _*$ and $F_*$ from closed subgroups
$\theta _*^{-1}(W^MN)_{\alpha , \beta }$ and $F_*^{-1}(W^ME)_{\alpha
, \beta }$ in $(W^MN_1)_{\alpha , \beta }$ and $(W^ME_1)_{\alpha ,
\beta }$ respectively onto $(W^MN)_{\alpha , \beta }$ and
$(W^ME)_{\alpha , \beta }$ by closed subgroups $\theta _*^{-1}(e)$
and $F_*^{-1}(e)$ correspondingly.

\par {\bf 3.} We have that $g \in C^{\alpha }_{\beta }({\bar M},
\{ s_{0,q}: q=1,...,k \};{\cal W},y_0)$ if and only if $f_2\circ g
\circ f_1^{-1}\in C^{\alpha }_{\beta }({\bar M}_1,\{ s_{0,q,1}:
q=1,...,k \} ; {\cal W}_1,y_1)$, where $f_1(s_{0,q})=s_{0,q,1}$ for
each $q=1,...,k$, $f_2(y_0\times e)=y_1\times e$. At the same time
$\psi \in Di^{\alpha }_{\beta }({\bar M})$ if and only if $f_1\circ
\psi \circ f_1^{-1}\in Di^{\alpha }_{\beta }({\bar M}_1)$ (see also
\S 3.2 \cite{lulal12}). Hence $(S^ME)_{\alpha , \beta }$ is
isomorphic with $(S^{M_1}E_1)_{\alpha , \beta }$ and inevitably wrap
groups $(W^ME)_{\alpha , \beta }$ and $(W^{M_1}E_1)_{\alpha , \beta
}$ are $C^{\alpha }_{\beta }$ diffeomorphic as differentiable spaces
and isomorphic as $C^l_{\beta }$ groups.

\par {\bf 12. Note.} Let $G$ be a topological group not necessarily associative, but
alternative:
\par $(A1)$ $g(gf)=(gg)f$ and $(fg)g=f(gg)$ and $g^{-1}(gf)=f$
and $(fg)g^{-1}=f$ for each $f, g\in G$ \\
and having a conjugation operation which is a continuous
automorphism of $G$ such that \par $(C1)$ $conj (gf)=conj (f)
conj(g)$ for each $g, f\in G$, \par $(C2)$ $conj (e)=e$ for the unit
element $e$ in $G$.
\par If $G$ is of definite class of smoothness, for example, $C^{\alpha }_{\beta }$
differentiable, then $conj$ is supposed to be of the same class. For
a commutative group in particular the identity mapping as the
conjugation can be taken. For $G= {\cal A}_r^*$ the usual
conjugation  $conj (z)={\tilde z}$ can be taken for each $z\in {\cal
A}_r^*$, where $1\le r\le 3$.
\par Suppose that \par $(A2)$ ${\hat G} =
{\hat G}_0u_0\oplus {\hat G}_1u_1\oplus ... \oplus {\hat
G}_{2^r-1}u_{2^r-1}$ such that $G$ is a multiplicative group of a
ring $\hat G$ with the multiplicative group structure, where ${\hat
G}_0,...,{\hat G}_{2^r-1}$ are pairwise isomorphic commutative
associative rings and $ \{ u_0,...,u_{2^r-1} \} $ are generators of
the Cayley-Dickson algebra ${\cal A}_r$ over a commutative field
$\bf K$, $1\le r\le 3$ and $(y_lu_l)(y_su_s)=(y_ly_s)(u_lu_s)$ is
the natural multiplication of any pure states in $G$ for $y_l\in
G_l$. For example, $G=({\cal A}_r^*)^n$ and ${\hat G} = {\cal
A}_r^n$.

\par {\bf 13. Lemma.} {\it If $G$ and $K$ are two topological
or differentiable groups twisted over $ \{ u_0,...,u_{2^r-1} \} $
satisfying conditions 12$(A1,A2,C1,C2)$ and $K$ is a closed normal
subgroup in $G$, where $2\le r\le 3$, then the quotient group is
topological or differentiable and twisted over $ \{
u_0,...,u_{2^r-1} \} $.}
\par {\bf Proof.} By the conditions of this lemma ${\hat G}= {\hat G}_0u_0\oplus {\hat G}_1u_1
\oplus ... \oplus {\hat G}_{2^r-1}u_{2^r-1}$, where ${\hat
G}_0,...,{\hat G}_{2^r-1}$ are pairwise isomorphic. Then ${\hat
G}/{\hat K} = ({\hat G}_0/{\hat K}_0)u_0\oplus ... \oplus ({\hat
G}_{2^r-1}/{\hat K}_{2^r-1})u_{2^r-1}$ is also twisted. Each ${\hat
G}_j$ is associative, hence $G/K$ is alternative, since $2\le r\le
3$ and using multiplicative properties of generators of the
Cayley-Dickson algebra ${\cal A}_r$. On the other hand, $conj
(K)=K$, hence $conj (gK)=K conj (g)= conj (g) K\in G/K$ and $conj
(ghK)= conj (gh) K=(conj (h) conj (g)) K = (conj (h)K) (conj (g)K) =
conj (hK) conj (gK)= conj (gK hK)$. \par The subgroup $K$ is closed
in $G$,
\par We recall the following definitions. If $G_1$ and $G_2$ are two differentiable groups, then their
product is supplied with the less fine plot structure for which
canonical projections $\pi _1: G\to G_1$ and $\pi _2: G\to G_2$ are
differentiable morphisms. That is a family $\sf P$ of plots of $G$
is such that $\pi _j\circ \sf P$ is an initial family ${\sf P}_j$ of
plots of $G_j$ for  $j=1$ and $j=2$. For differentiable groups as
usually $\sf P$ and ${\sf P}_j$ are preserved relative to the
inversions and multiplications in $G$ and $G_j$ respectively.
\par If $G$ and $F$ are two differentiable groups with families of
plots ${\sf P} = {\sf P}_G$ and ${\sf P}_F$ an algebraic morphism
$\theta : F\to G$ is a differentiable morphism if for each $h\in
{\sf P}_F$ the inclusion $\theta \circ h \in {\sf P}_G$ follows.
\par For two groups $G$ and $F$ with an algebraic embedding $\theta :
F\hookrightarrow G$ we supply $F$ with an induced differentiable
structure so that a family of plots ${\sf P}_F$ of $F$ is the less
fine for which $\theta $ is a differentiable morphism.
\par If $G$ is a differentiable group and $F$ is its algebraically
normal subgroup, we supply the quotient group $G/F$ with a
differentiable structure so that ${\sf P}_{G/F}$ is the most fine
plot structure for which the quotient mapping $\theta : G\to G/F$ is
a differentiable morphism.
\par Thus by the
definition of the quotient differentiable structure $G/K$ is the
differentiable group.

\par {\bf 14. Proposition.} {\it Let $\eta : N_1\to N_2$
be a $C^{\alpha '}_{\beta }$-retraction of $C^{\alpha '}_{\beta }$
differentiable spaces, $N_2\subset N_1$, $\eta |_{N_2}=id$, $y_0\in
N_2$, where $\alpha '\ge \alpha $, $\bar M$ is an $C^{\alpha
}_{\beta }$ differentiable space, $E(N_1,G,\pi ,\Psi )$ and
$E(N_2,G,\pi ,\Psi )$ are principal $C^{\alpha '}_{\beta }$ bundles
with a structure group $G$ satisfying conditions of \S \S 1, 2 in
Section 3 \cite{lulal12}. Then $\eta $ induces the group
homomorphism $\eta _*$ from $(W^ME;N_1,G,{\bf P})_{\alpha , \beta }$
onto $(W^ME;N_2,G,{\bf P})_{\alpha , \beta }$.}
\par {\bf Proof.} Due to Proposition 6$(1)$ the wrap group
$(W^ME;N_1,G,{\bf P})_{\alpha , \beta }$ is the principal $G^k$
bundle over $(W^MN_1)_{\alpha , \beta }$. We extend $\eta $ to
$\vartheta : E(N_1,G,\pi ,\Psi )\to E(N_2,G,\pi ,\Psi )$ such that
$\pi _2 \circ \vartheta =\eta \circ \pi _1$ and $pr_2\circ \vartheta
=id: G\to G$, where $pr_2: E_y\to G$ is the projection, $y\in N_1$,
$\eta (N_1)=N_2$. If $f\in C^{\alpha }_{\beta }({\bar M},N_1)$, then
$\eta \circ f := \eta (f(*))\in C^{\alpha }_{\beta }({\bar M},N_2)$.
If $f(s_{0,q})=y_0$, then $\eta (f(s_{0,q}))=y_0$, since $y_0\in
N_2$. From the inclusion $N_2\subset N_1$ we deduce that $C^{\alpha
}_{\beta }({\bar M},N_2)\subset C^{\alpha }_{\beta }({\bar M},N_1)$.
The parallel transport structure $\bf P$ is given over the same
differentiable space $\bar M$.
\par Put $\eta _* (<{\bf P}_{{\hat {\gamma }},u}>_{\alpha , \beta }) = <{\bf
P}_{\eta \circ {\hat {\gamma }},u}>_{\alpha , \beta }$, where ${\hat
{\gamma }}: {\hat M}\to N_1$. In accordance with Theorems 2.3 and
2.6 in Section 3 \cite{lulal12} $\eta _*(<{\bf P}_{{\hat {\gamma
}}_1,u}\vee {\bf P}_{{\hat {\gamma }}_2,u}>_{\alpha , \beta } = \eta
^*(<{\bf P}_{{\hat {\gamma }}_1,u}>_{\alpha , \beta }) \eta _*(<{\bf
P}_{{\hat {\gamma }}_2,u}>_{\alpha , \beta })$, and we can put $\eta
_*(q^{-1})= [\eta _*(q)]^{-1}$, consequently, $\eta _*$ is the group
homomorphism. Moreover, for each $g\in (W^ME;N_2,G,{\bf P})_{\alpha
, \beta }$ there exists $q\in (W^ME;N_1,G,{\bf P})_{\alpha , \beta
}$ such that $\eta _*(q)=g$, since $\gamma : {\bar M}\to N_2$ and
$N_2\subset N_1$ imply $\gamma : {\bar M}\to N_1$. On the other
hand, the structure group $G$ is the same, hence $\eta _*$ is the
epimorphism.

\par {\bf 15. Definition.} Suppose that $G$ is a topological group
satisfying Conditions 12$(A1,A2,C1,C2)$ such that $G$ is a
multiplicative group of the ring $\hat G$, where $1\le r\le 2$. We
define a smashed product $G^s$ such that it is a multiplicative
group of the ring ${\hat G}^s := {\hat G}\otimes _l{\hat G}$, where
$l=u_{2^r}$ denotes the doubling generator, a multiplication in
${\hat G}\otimes _l{\hat G}$ is given by the formula: \par $(1)$
$(a+bl)(c+vl) = (ac - v^* b) + (va+bc^*)l$ for each $a, b, c, v\in
\hat G$, where $v^* = conj (v)$.
\par In this relation it is worth to mention that ${\cal A}_r$ with
$r\le 3$ is the division algebra. For matrices with entries in
${\cal A}_r$ the Gauss' algorithm is valid, so matrices have ranks
by rows and columns which coincide and so a dimension over ${\cal
A}_r$ is defined \cite{dickselpap}.
\par A smashed product $M_1\otimes _lM_2$ of $C^{\alpha }_{\beta }$ differentiable manifolds $M_1, M_2$
over ${\cal A}_r$ with $dim_{{\cal A}_r} (M_1)= dim_{{\cal A}_r}
(M_2)$ is defined to be an ${\cal A}_{r+1}$ differentiable manifold
with local coordinates $z=(x,yl)$ of class $C^{\alpha }_{\beta }$,
where $x$ in $M_1$ and $y$ in $M_2$ are local coordinates.

\par Its existence and detailed description are demonstrated below.

\par {\bf 16. Proposition.} {\it The ring ${\hat G}^s$ from \S 15
has a multiplicative group $G^s$ containing all $a+bl\ne 0$ with $a,
b\in  {\hat G}$. If $\hat G$ is a topological or $C^{\alpha }_{\beta
}$ differentiable ring over ${\cal A}_r$, then ${\hat G}^s$ is a
topological or $C^{\alpha }_{\beta }$ differentiable over ${\cal
A}_{r+1}$ ring.}
\par {\bf Proof.} If $1\le r \le 2$ then a group $G$ is
associative, since the generators $ \{ u_0,...,u_{2^r-1} \} $ form
the associative group, when $r\le 2$. An element $a+bl\in {\hat
G}^s$ is non-zero if and only if $(a+bl)(a+bl)^* = aa^* +bb^*\ne 0$
due to 12$(A1,A2,C1,C2)$ and 15$(1)$. For $a+bl\ne 0$ we put $u =
(a^* - lb^*)/(aa^*+bb^*)$, where $aa^*+bb^*\in G_0$, hence $u (a+bl)
= (a+bl) u= 1\in G_0$, since $G_j$ is commutative for each
$j=0,...,2^r-1$, where $G_j$ denotes the multiplicative group of the
ring ${\hat G}_j$. The family of generators $ \{ u_0,...,
u_{2^{r+1}-1} \} $ for $r\le 2$ forms the alternative group, hence
${\hat G}^s = {\hat G}_0u_0\oplus ... \oplus {\hat
G}_{2^{r+1}-1}u_{2^{r+1}-1}$ is alternative, where ${\hat G}_j$ are
isomorphic with ${\hat G}_0$ for each $j$.
\par If an operation of the addition in $\hat G$ is continuous, then evidently
$(a+bl) + (c+ql)= (a+c) +(b+q)l$ is continuous. If an operation of
the multiplication in $\hat G$ is continuous, then Formula 15$(1)$
shows that the multiplication in ${\hat G}^s$ is continuous as well.
\par We have the decomposition ${\cal A}_{r+1} = {\cal A}_r\oplus
{\cal A}_rl$. If $\hat G$ is $C^{\alpha }_{\beta }$ differentiable,
then from the definition of plots it follows, that ${\hat G}^s$ is
$C^{\alpha }_{\beta }$ differentiable over ${\cal A}_{r+1}$ (see
also in details 17$(1-5)$).

\par {\bf 17. Theorem.} {\it Let ${\bar M}_1, {\bar M}_2$ and $N_1, N_2$ be
$C^{\alpha }_{\beta }$ differentiable manifolds over ${\cal A}_r$
with $1\le r \le 2$, and let $G$ be a group satisfying Conditions
12$(A1,A2,C1,C2)$. Suppose also that ${\bar M}_1\otimes _l{\bar
M}_2$, $N_1\otimes _lN_2$ are smashed products of differentiable
manifolds and $G^s$ is a smashed product group (see Proposition 16),
where $dim_{{\cal A}_r} ({\bar M}_1) = dim_{{\cal A}_r} ({\bar
M}_2)$, $dim_{{\cal A}_r} (N_1) = dim_{{\cal A}_r} (N_2)$. Then the
wrap group \par $(W^{M_1\otimes _lM_2; \{ s_{0,j,1}\otimes _l
s_{0,v,2}: j=1,...,k_1; v=1,...,k_2 \} } E;N_1\otimes _lN_2,G^s,{\bf
P}^s)_{\alpha , \beta }$ is twisted over $ \{ u_0,...,u_{2^{r+1}-1}
\} $ and is isomorphic with the smashed product
\par $(W^{M_2; \{ s_{0,v,2}: v=1,...,k_2 \} }E;N_1,(W^{M_1; \{
s_{0,j,1}: j=1,...,k_1 \} }E;N_1,G,{\bf P}_1)_{\alpha , \beta },{\bf
P}_2)_{\alpha , \beta }\otimes _l$
\par $(W^{M_2; \{ s_{0,v,2}: v=1,...,k_2 \} }E;N_2,(W^{M_1; \{
s_{0,j,1}: j=1,...,k_1 \} }E;N_2,G,{\bf P}_1)_{\alpha , \beta },{\bf
P}_2)_{\alpha , \beta }$ \\ of twice iterated wrap groups twisted
over $ \{ u_0,...,u_{2^r-1} \} $.}

\par {\bf Proof.} Let $M_b$ and $N_b$ be $C^{\alpha }_{\beta }$ differentiable manifolds over
${\cal A}_r$ with $1\le r \le 2$, $b=1, 2$ and let $G$ be a group
satisfying Conditions 12$(A1,A2,C1,C2)$ such that $E(N_b,G,\pi ,\Psi
)$ is a principal $G$-bundle. We consider the smashed products
${\bar M}_1\otimes _l{\bar M}_2$, $N_1\otimes _lN_2$ of $C^{\alpha
}_{\beta }$ differentiable manifolds and the smashed product group
$G^s$ (see Proposition 16). For $At ({\bar M}_b) = \{ (U_{j,b},\phi
_{j,b}): j \} $ an atlas of ${\bar M}_b$ its connecting mappings
$\phi _{j,b}\circ \phi _{k,b}^{-1}$ are $C^{\alpha }_{\beta }$
functions over ${\cal A}_r$ for $U_{j,b}\cap U_{k,b}\ne \emptyset $,
where $\phi _{j,b}: U_{j,b}\to {\cal A}_r$ are homeomorphisms of
$U_{j,b}$ onto $\phi _{j,b}(U_{j,b})$. Then ${\bar M}_1\otimes
_l{\bar M}_2$ consists of all points $(x,yl)$ with $x\in {\bar M}_1$
and $y\in {\bar M}_2$, with the atlas $At ({\bar M}_1\otimes _l{\bar
M}_2) = \{ (U_{j,1}\otimes _lU_{q,2}, \phi _{j,1}\otimes _l\phi
_{q,2}): j, q \} $ such that $\phi _{j,1}\otimes _l \phi _{q,2}:
U_{j,1}\otimes _lU_{q,2}\to {\cal A}_{r+1}^m$, where $m$ is a
dimension of ${\bar M}_1$ over ${\cal A}_r$. Express for $z=x+yl\in
{\cal A}_{r+1}$ with $x, y\in {\cal A}_r$ numbers $x, y$ in the $z$
representation. Then we denote by $\theta _{j,q}$ mappings
corresponding to $\phi _{j,1}\otimes _l\phi _{q,2}$ in the $z$
representation. Thus the transition mappings $\theta _{j,q}\circ
\theta _{k,n}^{-1}$ are $C^{\alpha }_{\beta }$ differentiable over
$\bf K$, when $(U_{j,1}\otimes _lU_{q,2})\cap (U_{k,1}\otimes
_lU_{n,2})\ne \emptyset $. Therefore, ${\bar M}_1\otimes _l{\bar
M}_2$ and $N_1\otimes _lN_2$ are $C^{\alpha }_{\beta }$
differentiable manifolds over ${\cal A}_{r+1}$.

\par Each $C^{\alpha }_{\beta }$ function $f(x,y) = f_1(x,y) + f_2(x,y)l$ by
$x\in U$ and $y\in V$ is of class $C^{\alpha }_{\beta }$ by
variables over $\bf K$ and takes values in $X_N\oplus X_Nl$ over
${\cal A}_{r+1}$, where $U$ and $V$ are open in $X_M$, $f_b(x,y)$ is
a $C^{\alpha }_{\beta }$ function with values in $X_N$ over ${\cal
A}_r$, $b=1, 2$.  If $z\in {\cal A}_{r+1}$, then
\par $(1)$ $z= v_0u_0+...+v_{2^{r+1}-1} u_{2^{r+1}-1}$, where
$v_j\in \bf K$ for each $j=0,...,2^{r+1}-1$, \par $(2)$ $v_0 = (z +
(2^{r+1}-2)^{-1} \{ -z + \sum_{j=1}^{2^{r+1}-1} u_j(zu_j^*) \} )/2$,
\par $(3)$ $v_s = (u_s(2^{r+1}-2)^{-1} \{ - z + \sum_{j=1}^{2^{r+1}-1}
u_j(zu_j^*) \} - z u_s )/2$ for each $s=1,...,2^{r+1}-1$, where $z^*
= {\tilde z}$ denotes the conjugated Cayley-Dickson number $z$. At
the same time we have for $z = x+yl$ with $x, y\in {\cal A}_r$, that
\par $(4)$ $x = v_0u_0+...+v_{2^r-1}u_{2^r-1}$ and \par $(5)$ $y =
(v_{2^r}u_{2^r}+...+v_{2^{r+1}-1}u_{2^{r+1}-1})l^*$, \\ where
$l=u_{2^r}$ denotes the doubling generator. Therefore, using
Formulas $(1-5)$ we get, that $f(x,y)$ is $C^{\alpha }_{\beta }$
differentiable over $\bf K$.

\par Then $E(N_1\otimes _lN_2,G^s,\pi ^s,\Psi ^s)$ is naturally isomorphic
with $E(N_1,G,\pi _1,\Psi _1)\otimes _l E(N_2,G,\pi _2,\Psi _2)$,
where $\pi ^s =\pi _1\otimes \pi _2l: E(N_1\otimes _lN_2,G^s,\pi
^s,\Psi ^s)\to N_1\otimes _l N_2$ is the natural projection.

\par If $\gamma : {\bar M}_1\otimes _l{\bar M}_2\to N_1\otimes _lN_2$ is an
$C^{\alpha }_{\beta }$ mapping, then $\gamma (z) = \gamma _1(x,y)
\times \gamma _2(x,y)l$, where $x\in {\bar M}_1$ and $y\in {\bar
M}_2$, $z=(x,yl)\in {\bar M}_1\otimes _l{\bar M}_2$, $\gamma _b:
{\bar M}_1\otimes _l{\bar M}_2\to N_b$. We can write $\gamma
_b(x,y)$ as $(\gamma _{b,1}(x))(y)$ a family of functions by $x$ and
a parameter $y$ or as $(\gamma _{b,2}(y))(x)$ a family of functions
by $y$ with a parameter $x$. If $\eta _{b,a}: {\bar M}_a\to N_b$,
then ${\bf P}_{{\hat {\eta }}_{b,a},u_b,a}$ denotes the parallel
transport structure on ${\bar M}_a$ over $E(N_b,G,\pi _b,\Psi _b)$.

\par Then we obtain \par ${\bf P}^s_{{\hat {\gamma }},u}
(z) = [{\bf P}_{{\hat {\gamma }}_{1,1},u_1;1}(x)] [{\bf P}_{{\hat
{\gamma }}_{1,2},u_1;2}(y)] \otimes _l [{\bf P}_{{\hat {\gamma
}}_{2,1},u_2;2}(x)] [{\bf P}_{{\hat {\gamma }}_{2,2},u_2;2}(y)] \in
E_{y_0}(N_1\otimes _lN_2,G^s,\pi ^s, \Psi ^s)$ \\ is the parallel
transport structure in ${\bar M}_1\otimes _l{\bar M}_2$ induced by
that of in ${\bar M}_1$ and ${\bar M}_2$, where $u\in E_{y_0}
(N_1\otimes _lN_2, G^s, \pi ^s, \Psi ^s)$, $u= u_1\otimes _l u_2$,
$u_b\in E_{y_{0,b}}(N_b,G,\pi _b,\Psi _b)$, $y_{0,b}\in N_b$ is a
marked point, $b=1, 2$, $y_0 = y_{0,1}\otimes _l y_{0,2}$. Thus
${\bf P}^s$ is $G^s$ equivariant. Therefore, the formula $<{\bf
P}^s_{{\hat {\gamma }},u}>_{\alpha , \beta } = < {\bf P}_{{\hat
{\gamma }}_1,u_1}>_{\alpha , \beta } \otimes _l <{\bf P}_{{\hat
{\gamma }}_2,u_2}>_{\alpha , \beta } = <[{\bf P}_{{\hat {\gamma
}}_{1,1},u_1;1}(x)] [{\bf P}_{{\hat {\gamma
}}_{1,2},u_1;2}(y)]>_{\alpha , \beta } \otimes _l <[{\bf P}_{{\hat
{\gamma }}_{2,1},u_2;2}(x)] [{\bf P}_{{\hat {\gamma
}}_{2,2},u_2;2}(y)]>_{\alpha , \beta }$, where ${\bf P}_{{\hat
{\gamma }}_b,u_b}$ produces the parallel transport structure in
${\bar M}_1\otimes _l{\bar M}_2$ over $E(N_b,G,\pi _b,\Psi _b)$,
$b=1, 2$.

\par Hence the group $(W^{M_1\otimes _lM_2; \{
s_{0,j,1}\otimes _l s_{0,v,2}: j=1,...,k_1; v=1,...,k_2 \} }
E;N_1\otimes _lN_2,G^s,{\bf P}^s)_{\alpha , \beta }$ is isomorphic
with the smashed product \par $(W^{M_2; \{ s_{0,v,2}: v=1,...,k_2 \}
}E;N_1,(W^{M_1; \{ s_{0,j,1}: j=1,...,k_1 \} }E;N_1,G,{\bf
P}_1)_{\alpha , \beta },{\bf P}_2)_{\alpha , \beta }\otimes _l$ \par
$ (W^{M_2; \{ s_{0,v,2}: v=1,...,k_2 \} }E;N_2,(W^{M_1; \{
s_{0,j,1}: j=1,...,k_1 \} }E;N_2,G,{\bf P}_1)_{\alpha , \beta },{\bf
P}_2)_{\alpha , \beta }$ \\ of iterated wrap groups.

\par {\bf 18. Theorem.} {\it A homomorphism of iterated
wrap groups $\theta : (W^ME)_{a;\infty ,\beta }\otimes
(W^ME)_{b;\infty ,\beta }\to (W^ME)_{a+b;\infty ,\beta }$ exists for
each $a, b\in \bf N$, where $G$ is a $C^{\infty }_{\beta }$ group,
$E(N,G,\pi ,\Psi )$ is a principal $C^{\infty }_{\beta }$ bundle
with a structure group $G$. Moreover, if $G$ is either associative
or alternative, then the homomorphism $\theta $ is either
associative or alternative correspondingly.}

\par {\bf Proof.} We consider iterated wrap groups $(W^ME)_{a;\infty ,\beta }$
as in \S 3, $a\in \bf N$. If $\gamma _a: {\bar M}^a\to N$, $\gamma
_b: {\bar M}^b\to N$ are $C^{\infty }_{\beta }$ mappings such that
$\gamma _b(s_{0,j_1}\times ...\times s_{0,j_b})=y_0$ for each
$j_l=1,...,k$ and $l=1,...,b$, then $\gamma := \gamma _a\times
\gamma _b: {\bar M}^a\times {\bar M}^b\to N\times N=N^2$, where
${\bar M}^a\times {\bar M}^b = {\bar M}^{a+b}$, $s_{0,j}$ are marked
points in $M$ with $j=1,...,k$ and $y_0$ is a marked point in $N$,
$C^{\infty }_{\beta } = \bigcap_{\alpha \in \bf N}C^{\alpha }_{\beta
}$. This gives the iterated parallel transport structure ${\bf
P}_{{\hat {\gamma }},u;a+b}(x) := {\bf P}_{{\hat {\gamma
}}_a,u_a;a}(x_a)\otimes {\bf P}_{{\hat {\gamma }}_b,u;b}(x_b)$ on
${\bar M}^{a+b}$ over $E(N^2,G^2,\pi ,\Psi )$, where $u_b\in
E_{y_0}(N,G,\pi ,\Psi )$, $u = u_a\times u_b \in E_{y_0\times y_0}
(N^2,G^2,\pi ,\Psi )$.
\par The bunch ${\bar M}^b\vee {\bar M}^b$ is taken by points $s_{j_1,...,j_b}$
in ${\bar M}^b$, where $s_{j_1,...,j_b} := s_{0,j_1}\times ...
\times s_{0,j_b}$ with $j_1,...,j_b\in \{ 1,...,k \} $; $s_{0,j}$
are marked points in $\bar M$ with $j=1,...,k$. Then the
differentiable space $({\bar M}^a\vee {\bar M}^a)\times ({\bar
M}^b\vee {\bar M}^b)\setminus \{ s_{j_1,...,j_{a+b}}: j_l=1,...,k;
l=1,...,a+b \} $ is $C^{\alpha }_{\beta }$ homeomorphic with ${\bar
M}^{a+b}\vee {\bar M}^{a+b}\setminus \{ s_{j_1,...,j_{a+b}}:
j_l=1,...,k; l=1,...,a+b \} $, since $s_{j_1,...,j_a}\times
s_{j_{a+1},...,j_{a+b}} = s_{j_1,...,j_{a+b}}$ for each
$j_1,...,j_{a+b}$. We have also the embedding $Di^{\infty }_{\beta
}({\bar M}^a)\times Di^{\infty }_{\beta }({\bar M}^b)\hookrightarrow
Di^{\infty }_{\beta }({\bar M}^{a+b})$ for each $a, b\in \bf N$ (see
also \S 3.2 \cite{lulal12}). If $f_a \in Di^{\infty }_{\beta }({\bar
M}^a)$ having a restriction $f_a|_{K_a}=id$, then $f_a\times f_b\in
Di^{\infty }_{\beta }({\bar M}^{a+b})$ and $f_a\times
f_b|_{K_a\times K_b}=id$ for $K_a\subset M^a$. Put \par $\theta (
<{\bf P}_{{\hat {\gamma }}_a,u_a;a}>_{K,\infty ,\beta ;a}, <{\bf
P}_{{\hat {\gamma }}_b,u_b;b}>_{K,\infty ,\beta ;b}) =$ \par $ <
<{\bf P}_{{\hat {\gamma }}_a,u_a;a}>_{K,\infty ,\beta ;a}\otimes
<{\bf P}_{{\hat {\gamma }}_b,u_b;b}>_{K, \infty ,\beta b}>_{K,\infty
,\beta ;a+b},$
\\ so it is the group homomorphism, where the detailed notation
$<*>_{K,\alpha ,\beta ;a}$ means the equivalence class over the
differentiable space ${\bar M}^a$ instead of $\bar M$, $a\in \bf N$.

\par Therefore,
$<{\bf P}_{{\hat {\gamma }}\vee {\hat {\eta }},u;a+b}>_{K,\infty
,\beta ;a+b} :=$ \par $ < <{\bf P}_{{\hat {\gamma }}_a\vee {\hat
{\eta }}_a,u_a;a}>_{K,\infty ,\beta ;a}\otimes <{\bf P}_{{\hat
{\gamma }}_b\vee
{\hat {\eta }}_b,u_b;b}>_{K,\infty ,\beta ;b}>_{K,\infty ,\beta ;a+b} $ \\
$=< (<{\bf P}_{{\hat {\gamma }}_a,u_a;a}>_{K,\infty ,\beta ;a} <{\bf
P}_{{\hat {\eta }}_a,u_a;a}>_{K,\infty ,\beta ;a})\otimes  (<{\bf
P}_{{\hat {\gamma }}_b,u_b;b}>_{K,\infty ,\beta ;b} <{\bf P}_{{\hat
{\eta
}}_b,u_b;b}>_{K,\infty ,\beta ;b})>_{K,\infty ,\beta ;a+b}$ \\
$=< (<{\bf P}_{{\hat {\gamma }}_a,u_a;a}>_{K,\infty ,\beta
;a}\otimes <{\bf P}_{{\hat {\gamma }}_b,u_b;b}>_{K,\infty ,\beta
;b}) (<{\bf P}_{{\hat {\eta }}_a,u_a;a}>_{K,\infty ,\beta ;a}\otimes
<{\bf P}_{{\hat {\eta
}}_b,u_b;b}>_{K,\infty ,\beta ;b})>_{K,\infty ,\beta ;a+b}$ \\
$=< <{\bf P}_{{\hat {\gamma }}_a,u_a;a}>_{K,\infty ,\beta ;a}\otimes
<{\bf P}_{{\hat {\gamma }}_b,u_b;b}>_{K,\infty ,\beta ;b}>_{K,\infty
,\beta ;a+b} < <{\bf P}_{{\hat {\eta }}_a,u_a;a}>_{K,\infty ,\beta
;a}\otimes <{\bf P}_{{\hat {\eta
}}_b,u_b;b}>_{K,\infty ,\beta ;b} >_{K,\infty ,\beta ;a+b}$ \\
$=\theta ( <{\bf P}_{{\hat {\gamma }}_a,u_a;a}>_{K,\infty ,\beta
;a}, <{\bf P}_{{\hat {\gamma }}_b,u_b;b}>_{K,\infty ,\beta ;b})
\theta ( <{\bf P}_{{\hat {\eta }}_a,u_a;a}>_{K,\infty ,\beta ;a},
<{\bf P}_{{\hat {\eta
}}_b,u_b;b}>_{K,\infty ,\beta ;b})$. \\
Thus $\theta $ is the group homomorphism.
\par The mapping $C^{\infty }_{\beta }({\bar M}^a,N)\times C^{\infty }_{\beta }({\bar M}^b,N)\ni
(\gamma _a\times \gamma _b)\mapsto (\gamma _a,\gamma _b)\in
C^{\infty }_{\beta }({\bar M}^{a+b},N^2)$ is of $C^{\infty }_{\beta
}$ class. The multiplication in the group $G^v$ is defined by the
formula: $(a_1,...,a_v)\times (b_1,...,b_v)=(a_1b_1,...,a_vb_v)$,
where $G^v$ is the $v$ times direct product of $G$, $a_1,...,a_v,
b_1,...,b_v\in G$. Therefore, the multiplication in $G^v$ is
$C^{\infty }_{\beta }$ smooth for each $v\in \bf N$, since it is
such in $G$.
\par The iterated wrap group $(W^ME)_{l;\alpha ,\beta }$ for the bundle $E$ is
the principal $G^{kl}$ bundle over the iterated commutative wrap
group $(W^MN)_{l;\alpha ,\beta }$ for the manifold $N$, since the
number of marked points in $M^l$ is $kl$, where $E$ is the principal
$G$ bundle on the manifold $N$, $l\in \bf N$. Thus the iterated wrap
group is associative or alternative if such is $G$. In view of
Proposition 6 and \S 3 the homomorphism $\theta $ is of $C^{\infty
}_{\beta }$ class. From the wrap monoids it has the natural
$C^{\infty }_{\beta }$ extension on wrap groups.
\par If $G$ is associative, then \par $<{\bf P}_{{\hat {\gamma
}},u;a+b+v}>_{K,\infty ,\beta ;a+b+v} = < < ( <{\bf P}_{{\hat
{\gamma }}_a,u_a; a}>_{K,\infty ,\beta ;a} \otimes <{\bf P}_{{\hat
{\gamma }}_b,u_b;b}>_{K,\infty ,\beta ;b})
>_{K,\infty ,\beta ;a+b} \otimes <{\bf P}_{{\hat {\gamma }}_v,u_v; v}
>_{K,\infty ,\beta ;v}>_{K,\infty ,\beta ;a+b+v}$ \\ $ = <
<{\bf P}_{{\hat {\gamma }}_a,u_a;a}>_{K,\infty ,\beta ;a}\otimes (
<{\bf P}_{{\hat {\gamma }}_b,u_b;b}>_{K,\infty ,\beta ;b}\otimes
<{\bf P}_{{\hat {\gamma }}_v,u_v; v}>_{K,\infty ,\beta ;v})
>_{K,\infty ,\beta ;a+b+v} = \theta ( \theta ( <{\bf P}_{{\hat {\gamma
}}_a,u_a; a}>_{K,\infty ,\beta ;a}, <{\bf P}_{{\hat {\gamma
}}_b,u_b;b}>_{K,\infty ,\beta ;b}), <{\bf P}_{{\hat {\gamma
}}_v,u_v; v}>_{K,\infty ,\beta ;v}) $  \\
$\theta ( <{\bf P}_{{\hat {\gamma }}_a,u_a; a}>_{K,\infty ,\beta
;a}, \theta ( <{\bf P}_{{\hat {\gamma }}_b,u_b;b}>_{K,\infty ,\beta
;b}), <{\bf P}_{{\hat
{\gamma }}_v,u_v; v}>_{K,\infty ,\beta ;v})) $, \\
consequently, $\theta $ is the associative homomorphism.
\par If $G$ is alternative, then
\par $<{\bf P}_{{\hat {\gamma }},u;a+a+b}>_{K,\infty ,\beta ;a+a+b} = < < (
<{\bf P}_{{\hat {\gamma }}_a,u_a; a}>_{K,\infty ,\beta ;a} \otimes
<{\bf P}_{{\hat {\gamma }}_a,u_a;a}>_{K,\infty ,\beta ;a})
>_{K,\infty ,\beta ;a+a} \otimes <{\bf P}_{{\hat {\gamma }}_b,u_b; b}
>_{K,\infty ,\beta ;v}>_{K,\infty ,\beta ;a+a+b}$ \\ $ = <
<{\bf P}_{{\hat {\gamma }}_a,u_a;a}>_{K,\infty ,\beta ;a}\otimes (
<{\bf P}_{{\hat {\gamma }}_a,u_a;a}>_{K,\infty ,\beta ;a}\otimes
<{\bf P}_{{\hat {\gamma }}_b,u_b; b}>_{K,\infty ,\beta ;b})
>_{K,\infty ,\beta ;a+a+b} = \theta ( \theta ( <{\bf P}_{{\hat {\gamma
}}_a,u_a; a}>_{K,\infty ,\beta ;a}, <{\bf P}_{{\hat {\gamma
}}_a,u_a;a}>_{K,\infty ,\beta ;a}), <{\bf P}_{{\hat {\gamma
}}_b,u_b; b}>_{K,\infty ,\beta ;b}) = $  \\
$\theta ( <{\bf P}_{{\hat {\gamma }}_a,u_a; a}>_{K,\infty ,\beta
;a}, \theta ( <{\bf P}_{{\hat {\gamma }}_a,u_a;a}>_{K,\infty ,\beta
;a}), <{\bf P}_{{\hat
{\gamma }}_b,u_b; b}>_{K,\infty ,\beta ;b})) $. \\
Thus the homomorphism $\theta $ is alternative from the left,
analogously it is alternative from the right.

\par {\bf 19. Remark.} Wrap groups were defined and studied above for
fiber bundles over a field $\bf K$ and algebras ${\cal A}_r$ with
$1\le r \le 3$ with a topological or differentiable structure group
$G$ which may be of Lie type as well.  In particular this
encompasses the case of multiplicative groups $G$ of commutative
algebras such as $H_n$ of diagonal matrices with entries in $\bf K$,
particularly, $H_4$ of quadra numbers. \par It is interesting to
mention that using an extension ${\bf F}\subset {\bf K}$ of a field
$\bf F$ up to a field $\bf K$ of a positive characteristic $p$ it is
possible to construct a division algebra of a dimension $m=n^2$
greater than eight over $\bf F$, where $n$ is a natural number
\cite{dickselpap}. The method of construction uses irreducible
polynomials over $\bf F$. This $\bf K$ may be a finite field
$F_{p^k}$ and a locally compact field of fractions $F_{p^k} (\theta
)$ over $F_{p^k}$ with an indeterminate $\theta $, where $k\in \bf
N$, that was outlined by Dickson. If such algebra is constructed
over $F_{p^k}$, then it exists over $F_{p^k}(\theta )$ also. Indeed,
each irreducible polynomial $g_l(x) := x^l + c_{l-1} x^{l-1} +...+
c_1x +c_0$ with expansion coefficients $c_0,...,c_{l-1}\in F_{p^k}$,
$c_0\ne 0$, is also irreducible over $F_{p^k}(\theta )$, since
$|c_j|=1$ for $c_j\ne 0$ and $|c_j|=0$ for $c_j=0$, where $|*|$ is
the multiplicative norm in $F_{p^k}(\theta )$ with $0<|\theta | <1$.
Each $x\in F_{p^k}(\theta )$ has the form $x=\theta ^N
\sum_{j=0}^{\infty } \beta _j \theta ^j$, where $\beta _j=\beta
_j(x)\in F_{p^k}$ for each $j$, $\beta _0\ne 0$, $N\in \bf Z$, $|x|=
|\theta |^N$. Therefore, $|g_l(x)| = |\theta |^{Nl}$ for $N<0$,
since $|a+b|\le \max (|a|, |b|)$ for all $a, b\in F_{p^k}(\theta )$.
For $N>0$ we get $|g_l(x)-c_0| = |\theta |^{Nk_0}$, where $k_0 =
\min \{ k: ~ c_k\ne 0, c_{k-1}=0,...,c_1=0 \} $, consequently,
$|g_l(x)|=|c_0|=1$ and such $x$ can not be a zero of $g_l(x)$. Thus
if the polynomial $g_l(x)$ has a zero $z\in F_{p^k}(\theta )$, then
$|z| =1$, since $c_0\ne 0$. Consider terms in an expansion of
$g_l(z)$ by degrees of $\theta $. For $\theta ^0$ the corresponding
coefficient $\beta _0=\beta _0(g_l(z))$ should be equal to zero,
hence a zero $x_0\in F_{p^k}$ would exist. Thus the polynomial
$g_l(x)$ is irreducible over $F_{p^k}(\theta )$ also.
\par Suppose now that $\sf R$ is a ring having two subgroups. One of them $G_1$ is
commutative related with the addition, $G_1=({\sf R},+)$. Another is
multiplicative $G_2= ({\sf R}\setminus \{ 0 \} ,\times )$.
Particularly, $\sf R$ may be an algebra $\sf A$ over the field $\bf
K$. We consider the cases of commutative, associative and as well as
non-associative rings and algebras with associative addition in
$G_1$ and alternative multiplication in $G_2$. Suppose that a fiber
bundle is given with the structure ring $\sf R$ or the structure
algebra $\sf A$ instead of a group. We suppose also that with $G_1$
and $G_2$ parallel transport structures $\mbox{}_1{\bf P}$ and
$\mbox{}_2{\bf P}$ are related. So we shall say that there is the
parallel transport structure $\bf P$ on the principal fiber bundle
$E(N,{\sf R},\pi , \Psi )$ or $E(N,{\sf A},\pi , \Psi )$.

\par {\bf 20. Theorem.} {\it Wrap groups $(W^{M, \{ s_{0,q}: q=1,...,k \} } E;N,{\sf R},
{\bf P})_{\alpha , \beta }$ or \\ $(W^{M, \{ s_{0,q}: q=1,...,k \} }
E;N,{\sf A},{\bf P})_{\alpha , \beta }$ exist with two $C^l_{\beta
}$ group's operations and they are the principal fiber bundles over
the commutative group $(W^MN)_{\alpha , \beta }$ with the structure
ring ${\sf R}^k$ or the structure algebra ${\sf A}^k$ respectively,
where $l=\alpha ' - \alpha $ for $\alpha \le \alpha '<\infty $ or
$l=\infty $ for $\alpha ' =\infty $ (see also Remark 19 above).}
\par {\bf Proof.} Earlier wrap groups
$(W^{M, \{ s_{0,q}: q=1,...,k \} } E;N,G_j,\mbox{}_j{\bf P})_{\alpha
, \beta }$ for $j=1, 2$ for the principal fiber bundles $E(N,G_j,\pi
_j, \Psi _j)$ having $C^l_{\beta }$ group's operations were
constructed in accordance with Theorem 6 \cite{lulal12}. We consider
$\theta _j: G_j\hookrightarrow {\sf R}$ or $\theta _j:
G_j\hookrightarrow {\sf A}$ group embedding of the $C^{\alpha
'}_{\beta }$ class of differentiability, $j=1, 2$. In view of
Proposition 6 above they are the principal fiber bundles over
$(W^MN)_{\alpha , \beta }$ with the structure groups $G_j^k$. At the
same time the principal fiber bundle $E(N,{\sf R},\pi , \Psi )$ or
$E(N,{\sf A},\pi , \Psi )$ is isomorphic with $E((E(N,G_2,\pi _2,
\Psi _2)),G_1,\pi _1, \Psi _1)/ \xi$, where the equivalence relation
$\xi $ is induced by the equality $\theta _1(x)=\theta _2(y)$ in
$\sf R$ or $\sf A$ respectively of the corresponding elements $x\in
G_1$ and $y\in G_2$. Then we put \par $(W^{M, \{ s_{0,q}: q=1,...,k
\} } E;N,{\sf T},{\bf P})_{\alpha , \beta } := $
\par $(W^{M, \{ s_{0,q}: q=1,...,k \} } E((E(N,G_2,\pi _2, \Psi
_2)),G_1,\pi _1, \Psi _1);{\bf P})_{\alpha , \beta }/\xi $, \\ where
${\sf T}={\sf R}$ or ${\sf T}={\sf A}$ respectively. In addition we
have put \par $(W^{M, \{ s_{0,q}: q=1,...,k \} }
E;K,G_j,\mbox{}_j{\bf P})_{\alpha , \beta } =: (W^{M, \{ s_{0,q}:
q=1,...,k \} } E(K,G_j,\pi _j, \Psi _j));{\bf P})_{\alpha , \beta }$, \\
consequently, $(W^{M, \{ s_{0,q}: q=1,...,k \} } E;N,{\sf T},{\bf
P})_{\alpha , \beta }$ is supplied with two $C^l_{\beta }$ group
operations corresponding to $(W^{M, \{ s_{0,q}: q=1,...,k \} }
E;N,G_j,\mbox{}_j{\bf P})_{\alpha , \beta }$ for $j=1, 2$ and there
exists the principal fiber bundle
\par $(W^{M, \{ s_{0,q}: q=1,...,k \} } E((E(N,G_2,\pi _2, \Psi _2)),G_1,
\pi _1, \Psi _1);{\bf P})_{\alpha , \beta }$
\par $\to (W^{M, \{ s_{0,q}: q=1,...,k
\} } E;N,G_2,\mbox{}_2{\bf P})_{\alpha , \beta }$ \\ with the
structure group $G_1^k$. Using the equivalence relation $\xi $
inevitably infers that there exists the principal fiber bundle \par
$(W^{M, \{ s_{0,q}: q=1,...,k \} } E;N,{\sf T},{\bf P})_{\alpha ,
\beta }\to (W^MN)_{\alpha , \beta }$
\\ with the structure ring ${\sf R}^k$ or the structure algebra
${\sf A}^k$ correspondingly.

\par {\bf 21. Remark.} Wrap groups and semigroups can be generalized
for an infinite discrete closed subset $ \{ s_{0,q}: q\in \lambda \}
$ of marked points in $\bar M$, when $\bar M$ is not compact, where
$\lambda $ is a set, $card (\lambda )\ge \aleph _0$.
\par Apart from the differentiable spaces $M$, $N$, $E$ considered
above over the infinite non-discrete field $\bf K$ analogous wrap
groups also exist when $M$, $N$ and $E$ are over a finite field
$F_{p^k}$ or an algebra ${\cal A}_r$ over $F_{p^k}$, $1\le r\le 3$,
but such wrap groups become already discrete. Particularly, if $M$,
$N$ and $E$ are finite, then wrap groups are finite (see also
References [25,26] in \cite{lulal12}). Therefore, we have considered
above topological infinite groups, when $\bf K$ is infinite and
non-discrete.

\par Ludkovsky S.V. Department of Applied Mathematics MIREA,
av. Vernadsky 78, Moscow 119454
\par sludkowski@mail.ru

\end{document}